\def\mapright#1{\smash{\mathop{\longrightarrow}\limits^{#1}}}

\documentclass{amsart}
\usepackage{amscd, amsmath, amsthm, amssymb}

\newtheorem{theorem}{Theorem}[section]
\newtheorem{lemma}[theorem]{Lemma}
\newtheorem{proposition}[theorem]{Proposition}

\newtheorem{corollary}[theorem]{Corollary}

\theoremstyle{definition}     % italic or bold etc.
\newtheorem{definition}[theorem]{Definition}

\newtheorem{example}[theorem]{Example}

\newtheorem{claim}[theorem]{Claim}
\newtheorem{question}[theorem]{Question}
\newtheorem{problem}[theorem]{Problem}
\theoremstyle{remark}
\newtheorem{remark}[theorem]{Remark}

\numberwithin{equation}{section}

\begin{document}

\title[hyperk\"aher manifolds ]
{Bimeromorphic automorphism groups of non-projective hyperk\"ahler manifolds  
- a note inspired by C. T. McMullen}

\author[K. Oguiso]{Keiji Oguiso}
\address{Graduate School of Mathematical Sciences,
University of Tokyo, Komaba, Meguro-ku,
Tokyo 153-8914, Japan
}
\email{oguiso@ms.u-tokyo.ac.jp}


\begin{abstract} 
Being inspired by a work of Curtis T. McMullen 
about a very impressive automorphism of 
a K3 surface of Picard 
number zero, we shall clarify the structure 
of the bimeromorphic automorphism group 
of a non-projective hyperk\"ahler manifold, up to 
finite group factor. We also discuss relevant topics, 
especially, new counterexamples of Kodaira's problem 
about algebraic approximation of a compact K\"ahler manifold. 
\end{abstract}

\maketitle

%%%%%%%%%%%%%%%%%%%%%%%%%%%%%%%%%%%%%%%%%%%

\setcounter{section}{0}
\section{Introduction - Background and the statement of main result} 
In his study about complex dynamics [Mc], Curtis T. McMullen 
has found a very impressive automorphism of 
a K3 surface of Picard 
number zero. In this note, inspired by his K3 automorphisms, 
we shall clarify the structure 
of the bimeromorphic automorphism group 
of a non-projective hyperk\"ahler manifold, up to finite group 
factor. We then discuss relevant topics, especially, 
new counterexamples of Kodaira's problem about algebraic approximation of 
a compact K\"ahler manifold. Our main results are 
Theorems (1.5), (1.9), (2.4) and (3.4). 
Throughout this note, we work in the category of complex varieties with 
Euclidean topology.  
\par \vskip 1pc

For the statement of McMullen's K3 automorphism (Theorem (1.1)), we first 
recall a few notions from his paper [Mc]. 
Two complex numbers $\alpha$ 
and $\beta$ are {\it multiplicatively independent} 
if the only solution to $ \alpha^{m}\beta^{n} = 1$ with $(m, n) 
\in \mathbf Z^{2}$ is $(0,0)$. 
An automorphism $F^{*}(z_{1}, z_{2}) = (\alpha_{1} z_{1}, \alpha_{2} z_{2})$ 
of the $2$-dimensional 
unit disk $\Delta^{2} := \{(z_{1}, z_{2}) \in \mathbf C^{2}\,  
\vert\, \vert z_{1} \vert < 1, \vert z_{2} \vert < 1\}$ is 
an {\it irrational rotation} if $\alpha_{1}$ and $\alpha_{2}$ are 
multiplicatively independent numbers on 
the unit circle $S^{1} := \{z \in \mathbf C\, \vert\, \vert z \vert = 1\}$. 
An $f$-stable domain $D \subset S$ is a {\it Siegel disk} of $f$ centered at 
$P$, if there are an analytic isomorphism 
$\varphi : (P, D) \simeq ((0,0), \Delta^{2})$ and an irrational rotation 
$F$ of $\Delta^{2}$ such that $f \vert D = 
\varphi^{-1} \circ F \circ \varphi$. As it is observed in [ibid], 
K3 surfaces having Siegel disk are 
never projective, or more strongly, of algebraic dimension $0$ (cf. Example (2.5)). 
\par \vskip 1pc

McMullen has found the following:

\begin{theorem} \label{theorem:mcmullen} {\rm [Mc]} There is a pair $(S, f)$ of a K3 
surface $S$ 
with $\rho(S) = 0$ 
and an automorphism $f$ of $S$ having a Siegel disk. In such a pair 
$(S, f)$, the topological entropy of $f$, that is, the natural logarithm of 
the 
spectral radius of 
$f^{*} \vert H^{2}(S, \mathbf Z)$, is always positive. 
Moreover, such pairs are at most countable. 
\end{theorem}  

\begin{definition} \label{definition:mcpair} {\it McMullen's pair} is a pair 
$(S, f)$ 
in Theorem (1.1). 
\end{definition}  

He constructed such a K3 automorphism starting from certain Salem polynomial 
of degree $22$. (See (3.1) for the definition of Salem polynomial.) 
Existence of McMullen's pair suggests some complexity of 
automorphisms of a non-projective K3 surface, 
and naturally leads us to the 
following: 

\begin{question} \label{question:nat} How complicated is the full 
automorphism group $\text{Aut}\, (S)$ of a K3 
surface $S$ with Picard number $\rho(S) = 0$? 
\end{question} 

The aim of this note 
is to study this question in a bit more general setting: 
"how does the bimeromorphic automorphism group ${\rm Bir}\, (M)$ of a non-projective hyperk\"ahler manifold 
$M$ look like?" This question is asked by Y. Kawamata to me (after a preliminary version of this note). 

\begin{definition} \label{definition:hk} {\it A hyperk\"ahler manifold} 
is a compact complex simply-connected K\"ahler 
manifold $M$ admitting an everywhere non-degenerate global 
holomorphic $2$-form $\sigma_{M}$ with $H^{0}(M, \Omega_{M}^{2}) = 
\mathbf C \sigma_{M}$. 
\end{definition} 

Hyperk\"ahler manifolds are even dimensional. They coincide with 
K3 surfaces in dimension $2$ by 
a result of Siu [Si], and share 
many of properties with K3 surfaces [Hu]. 
For instance, given a hyperk\"ahler manifold, its generic deformation 
is always non-projective, or more strongly, of trivial 
N\'eron-Severi group $\{0\}$ (see eg. [Og1, Corollary 1.3]). 

One of the most important properties of a hyperk\"ahler manifold is the existence of Beauville-Bogomolov-Fujiki's form 
(BF-form for short) $(*, **)$. BF-form is an integral symmetric bilinear from on 
$H^{2}(M, \mathbf Z)$ of signature $(3, 0, b_{2}(M) -3)$. (For more detail, 
see [Be1], [Fu3], an excellent survey [GHJ, Part III] by Huybrechts, and also Section 2.)  

By means of BF-form, the signature (i.e. the numbers of positive-, zero-, negative-eigenvalues) of 
the N\'eron-Severi group 
$NS(M)$ is either 
$$(1, 0, \rho(M) -1)\,\, ,\,\, (0, 1, \rho(M) -1)\,\, ,\,\, (0, 0, \rho(M))\,\, .$$ 
Here $\rho(M)$ is the Picard number of $M$. We call these three cases {\it hyperbolic}, {\it parabolic}, {\it elliptic} 
respectively. Due to a very deep result by Huybrechts [Hu], $M$ is projective iff $NS(M)$ is hyperbolic. 
 
Our main result is as follows: 

\begin{theorem} \label{theorem:main} The bimeromorphic automorphism group 
${\rm Bir}\,(M)$ of a non-projective hyperk\"ahler manifold is an almost abelian group of rank at most 
${\rm max}\,(\rho(M) -1, 1)$. It is, in particular, finitely generated. More precisely:

\begin{list}{}{
\setlength{\leftmargin}{10pt}
\setlength{\labelwidth}{6pt}
}
\item[(1)] If $NS(M)$ is elliptic, then ${\rm Bir}\,(M)$ 
falls into the exact sequence 
$$1 \longrightarrow N \longrightarrow {\rm Bir}\,(M) \longrightarrow 
\mathbf Z^{r} \longrightarrow 0$$ 
where $N$ is a finite group and $r$ is either $0$ or $1$.  Moreover, 
$r = 0$ if $b_{2}(M) - \rho(M)$ is odd, i.e. the rank of transcendental lattice (cf. Section 2) is odd.  
Moreover, in each dimension $2m$, both $r = 0$ and $r = 1$ are realizable. 
\item[(2)] If $NS(M)$ is parabolic, then ${\rm Bir}\,(M)$ is an almost abelian group 
of rank at most $\rho(M) -1$.  Moreover, in dimension $2$, this estimate is optimal. 
\end{list}
\end{theorem} 
Here an almost abelian group of rank $r$ is a group isomorphic to $\mathbf Z^{r}$ 
"up to finite kernel and cokernel" (See Section 9 for the precise 
definition employed here and some properties we need.) We also note that there 
is a projective hyperk\"ahler manifold $M$ s.t. ${\rm Bir}\, (M)$ is not almost abelian 
(cf. [Og2]). 
 
As a special case of Theorem (1.5), we answer Question (1.3) fairly completely:
\begin{corollary} \label{corollary:mc} Let $M$ be a K3 surface 
of Picard number $0$. Then ${\rm Aut}\,(M)$ is isomorphic to either 
$\{id\}$ or $\mathbf Z$. In particular, ${\rm Aut}\, (M) \simeq \mathbf Z$ for 
a McMullen's K3 surface $M$. 
Moreover, ${\rm Aut}\,(M) \simeq \mathbf Z$ for at most countably many 
$M$; in other words, mostly 
${\rm Aut}\,(M) = \{id\}$. 
\end{corollary} 
So, the automorphism group of a McMullen's K3 surface turns out to be quite simple 
as a group. It will be also very interesting to view these results from the view of topological 
entropy, especially, in connection with recent important results of Dinh and Sibony [DS, 1,2,3]. 
We also note that there is a projective hyperk\"ahler manifold $M$ whose ${\rm Bir}\, (M)$ is not 
almost abelian (cf. [Og2]). 
\par \vskip 1pc

%Key step is the finite generation of $\text{Aut}\,(M)$. 
%First we show that $\vert \text{Ker}\, (\chi) \vert < \infty$ by using a 
%celebrated solution of 
%Calabi's conjecture due to Shing-Tung Yau [Ya1] (see also [Ya2]) in 
%Theorem (2.4). 
%We then show, in Proposition (4.3),  
%the finite 
%generation of $\text{Im}\, (\chi)$  
%by using Dirichlet's unit theorem, one of the most fundamental theorems in 
%algebraic number theory (See eg. [Ta, Section 9.2, Theorem]). Finally, 
%in Proposition (4.4), we restrict 
%the size of 
%$\text{Im}\, (\chi)$ as claimed,  by 
%reducing 
%the problem to the finiteness result about Salem polynomials with bounded 
%degree and trace; see Definition (3.1) for the definition of a Salem 
%polynomial 
%and Theorem (3.3) and 
%Proposition (3.14) for necessary preparations. See also Remark (2.5) 
%about a few concrete examples of subgroups of $\mathbf C^{\times}$ 
%which should 
%be ruled out from the candidates of $\text{Im}\,(\chi)$. Though it will not 
%appear in this note, it would be very interesting to seek a relation 
%hidden between hyperk\"ahler manifolds and the finiteness of class numbers, 
%another most fundamental theorem in algebraic number theory.  
%(see [HLOY] about some relations with K3 surfaces.)  

As it is well known, there are hyperk\"ahler manifolds of dimension 
$2m$ for each $m \geq 2$; the Duady space $S^{[m]}$ of $0$-dimensional 
analytic 
subspaces $(Z, \mathcal O_{Z})$ of length $m$ on a K3 surface $S$ and 
manifolds $D_{m}$ obtained by a small deformation of $S^{[m]}$ ([Be1]). 
We have $b_{2}(D_{m}) = b_{2}(S^{[m]}) = 23$ by [ibid], 
and $\rho(D_{m}) = 0$ and $b_{2}(M) - \rho(M) = 23$, 
an odd number, for a 
sufficiently generic 
$D_{m}$ in the Kuranishi family (see eg. [Og]). As an immediate consequence 
of Theorem (1.5)(1), 
one has: 
\begin{corollary} \label{corollary:series} 
${\rm Bir}\,(D_{m})$ is a finite group for $D_{m}$ with $\rho(D_{m}) = 0$.  
\end{corollary} 

Compare Corollary (1.7) with (1.6). 
\par \vskip 1pc  

Next, we apply McMullen's pair to construct a simply-connected 
$d$-dimensional counterexample of Kodaira's problem about algebraic 
approximation of a compact K\"ahler manifold in each dimension $d \geq 4$: 

\begin{problem} \label{problem:cp} (cf. [CP]) Is any compact K\"ahler manifold 
algebraically 
approximated? Or more precisely, for a given compact K\"ahler manifold $X$, 
is there 
a small deformation $\pi : \mathcal X \longrightarrow \mathcal B$ 
such that $\mathcal X_{0} \simeq X$ and $\mathcal X_{t_{\nu}}$ is 
projective for a sequence $(t_{\nu})$ converging to $0$? 
\end{problem} 

\begin{theorem} \label{theorem:kd} Let $d \geq 4$ be any integer. Then there 
is a $d$-dimensional simply-connected compact 
K\"ahler manifold $Z_{d}$ which can not be algebraically approximated. 
Moreover, in our construction, $Z_{d}$ is bimeromorphic to a complex 
manifold 
having trivial 
canonical bundle if $d \not= 5$ and, in addition, is not analytically 
rigid if $d \not= 4, 5$. 
\end{theorem}

Theorem (1.9) is inspired by a very important work of C. Voisin 
[Vo] and a private communication with A. Fujiki 
about her work. She constructed 
counterexamples which are bimeromorphic to complex tori
in each dimension $d \geq 4$ and also simply-connected ones in each even 
dimension $2m \geq 6$. Her 
examples are, on the one hand, much stronger   
than just being counterexamples (see [ibid]). But,  
on the other hand, they are not simply-connected in odd dimension 
and are analytically rigid in even dimension. Our examples are made from a McMullen's pair, 
which {\it apriori} 
seems nothing to do with Kodaira's problem. 
\par \vskip 1pc 

%In Sections 2, after recalling the definition and basic properties 
%of BF-form, we study 
%two representations of the automorphism group 
%of a hyperk\"ahler manifold (Theorems (2.1) and (2.4)). In Section 3, 
%we clarify certain close relation 
%between automorphisms of a Picard generic hyperk\"ahler manifold and Salem 
%polynomials (Theorem (3.3); see also (3.1) for the definition of 
%Salem polynomial). In McMullen's pair $(S, f)$, the positivity of the 
%topological entropy of $f$ naturally leads us to the notion of Salem 
%polynomial. Indeed, McMullen constructed his pair starting from certain 
%Salem polynomial 
%of degree $22$ [Mc] (see also [BM]). This section is a sort of converse of 
%his construction in 
%hyperk\"ahler case. In Sections 4 we shall prove Theorem (1.5). We use the 
%results in Sections 2 and 3 in the proof. We shall then give two applications 
%for K3 surfaces in Section 5; one is about the finite generation of 
%$\text{Aut}\, (S)$ of a complex K3 surface 
%$S$ (Theorem (5.1)) and the other is about 
%a new behaviour of the automorphism group of a K3 surface under one 
%dimensional small deformation (Theorem (5.5)). In Section 6, we shall prove 
%Theorem (1.8) 
%by using Theorems (1.1) and (2.1).  

Section 2 is a preliminary section about some basic properties 
about hyperk\"ahler manifolds and their bimeromorphic automorphisms. We clarify a close relation 
between bimeromorphic automorphism in elliptic $NS(M)$ case and Salem polynomials in Section 3. 
We prove the main part of Theorem (1.5)(1) in Section 4, the main part of Theorem (1.5)(2) in Section 5. 
In Sections 6 and 7, we show Corollary (1.6), existence part and optimality 
part in dimension $2$ in Theorem (1.5). We prove Theorem (1.9) in Section 8. Section 9 
is a sort of appendix about almost abelian groups. 
\par \vskip 1pc

{\it Acknowledgement.} I would like to express my deep thanks 
to Professors A. Fujiki, D. Huybrechts, Y. Kawamata, J.H. Keum for their 
very valuable discussions and comments. I would like to express my thanks 
to Professors T.-C. Dinh, Johns H. Coates, I. Enoki, R. Goto, K. Ohno, T. Peternell 
for their interest in this work. The final version has been completed during my stay 
at KIAS in March 2005. I would like to express my thanks to Professors J. Hwang and B. Kim 
for invitation.

\par \vskip 1pc

\section{Some basic properties of Hyperk\"ahler manifolds}

In this preliminary section, first we recall the definition and some basics 
of Beauville-Bogomolov-Fujiki's form (BF-form) and bimeromorphic automorphisms from [Be1], [Fu3] 
and [GHJ, Part III]. We then study various representations of the bimeromorphic 
automorphism group of a hyperk\"ahler manifold. Theorem (2.3) is new.

{\bf 1.} Throughout Section 2, let $M$ be a hyperk\"ahler manifold of complex dimension $2m$ with 
$H^{0}(M, \Omega_{M}^{2}) 
= \mathbf C \sigma_{M}$. We may normailize $\sigma_{M}$ as 
$\int_{M} (\sigma_{M}\overline{\sigma}_{M})^{m} = 1$. We note that $H^{2}(M, \mathbf Z)$ 
is torsion free, because $\pi_{1}(M) = \{1\}$. Define a (real valued) 
quadratic form 
$\tilde{q}_{M}$ on $H^{2}(M, \mathbf Z)$ by
$$\tilde{q}_{M}(\alpha) := \frac{m}{2}\int_{M} 
\alpha^{2}(\sigma_{M}\overline{\sigma}_{M})^{m-1} + (1-m)(\int_{M} 
\alpha \sigma_{M}^{m-1}\overline{\sigma}_{M}^{m}) \cdot (\int_{M} 
\alpha\sigma_{M}^{m}\overline{\sigma}_{M}^{m-1})\,\, .$$ 
By [Be1] [Fu3], there is a unique positive constant $c_{M} > 0$ 
such that the symmetric bilinear form 
$$(*, **) : H^{2}(M, \mathbf Z) \times H^{2}(M, \mathbf Z) 
\longrightarrow \mathbf Z$$ 
associated with 
$q_{M} := c_{M}\tilde{q}_{M}$, is integral and primitive, in 
the sense that $(*, **)/k$ ($k \in \mathbf Z_{>0}$) is integral valued iff 
$k = 1$. {\it Beauville-Bogomolov-Fujiki's form}, BF-form for short, is 
this bilinear form 
$(*,**)$. BF-form $(*,**)$ coincides with the cup product when 
$\text{dim}\, M = 2$, i.e. when $M$ is a K3 surface, and 
$(*, **)$ 
and enjoys similar properties to the cup product of K3 surface.
\par \vskip 1pc

%By definition, $(*,**)$ is $\text{Aut}\, (M)$-stable in the sense that 
%$(f^{*}x, f^{*}y) = (x, y)$ for $\forall f \in \text{Aut}\,(M)$. Though less trivial, 
%BF-form is also stable under the bimeromorphic automorphism group ${\rm Bir}\, (M)$:

{\bf 2.} We always regard $H^{2}(M, \mathbf Z)$ as a lattice by 
BF-form. Then $H^{2}(M, \mathbf Z)$ is of signature $(3, 0, b_{2}(M) -3)$, 
i.e. the numbers of positive-, zero- and negative-eigenvalues are $3$, $0$, 
$b_{2}(M) - 3$. More precisely, 
by [Be1] one has that: 
$$(\sigma_{M}, \sigma_{M}) = 0\,\, ,\,\, 
(\sigma_{M}, \overline{\sigma}_{M}) > 0\,\, ,$$ 
$$(\sigma_{M}, H^{1,1}(M, \mathbf R)) = (\overline{\sigma}_{M}, 
H^{1,1}(M, \mathbf R)) = 0\,\, ,$$ 
and that if $\eta \in H^{1,1}(M, \mathbf R)$ is a K\"ahler class, then 
$(\eta, \eta) > 0$ and $(*,**)$ is negative definite on 
the orthogonal complement of $\mathbf R\langle 
\eta \rangle$ in $H^{1, 1}(M, \mathbf R)$. Thus BF-form on $H^{1,1}(M, \mathbf R)$ 
is of signature $(1, 0, h^{1,1}(M) -1)$, and $\{x \in H^{1,1}(M, \mathbf R) \vert (x^{2}) > 0\}$ 
has two connected components. We define the {\it positive cone} $\mathcal C = \mathcal C(M)$ to be the component which 
contains the K\"ahler cone $\mathcal K(M)$, i.e. the (open convex) cone of $H^{1,1}(M, \mathbf R)$ 
consisting of K\"ahler classes of $M$.

{\bf 3.} Let $NS(M)$ be the N\'eron-Severi group of $M$, i.e. the subgroup of $H^{2}(M, \mathbf Z)$ 
generated by the first Chern classes of holomorphic line bundles. We regard $NS(M)$ as a 
(possibly degenerate) sublattice of $H^{2}(M, \mathbf Z)$ by $(*, **)$. The {\it transcendental lattice} $T(M)$ of $M$ is the pair of the minimum  
subgroup $T$ of $H^{2}(M, \mathbf Z)$ such that $\mathbf C \sigma_{M} \in T \otimes \mathbf C$ and 
$H^{2}(M, \mathbf Z)/T$ is torsion-free, and the bilinear 
form $(*,**)_{T \times T}$. 
By the Lefschetz (1,1)-theorem and non-degeneracy of BF-form on $H^{2}(M, \mathbf Z)$, we have that 
$$NS(M) = H^{2}(M, \mathbf Z) \cap H^{1, 1}(M)\,\, ,\,\, 
T(M) = \{x \in H^{2}(M, \mathbf Z)\, \vert\, (x, NS(M)) = 0\}\, .$$ 
Note that $T(M) \otimes \mathbf R$ contains 
{\it a positive $2$-plane} 
$P := \mathbf R \langle \text{Re}\, \sigma_{M}\, , \text{Im}\, 
\sigma_{M}\rangle$. 
Then, as in the K3 case [Ni], the signature of $NS(M)$ is either one of 
$$(1\, ,\, 0\, ,\, \rho(M) -1)\,\, ,\,\,(0\, ,\, 1\, ,\, \rho(M) -1)\,\, 
,\,\, (0\, ,\, 0\, ,\, \rho(M))\,\, .$$ 
We call these three cases, {\it hyperbolic}, {\it parabolic}, {\it elliptic}, respectively. 
By a fundamental result of Huybrechts [Hu] (see also [GHJ, Proposition 26.13]), $M$ is projective iff 
$NS(M)$ is hyperbolic.

In the hyperbolic and elliptic cases, one has also
$$NS(M) \cap T(M) = \{0\}\,\, \text{and}\,\, [H^{2}(M, \mathbf Z) : 
NS(M) \oplus T(M)] < \infty\,\, .$$

In the parabolic case, $NS(M) \cap T(M) = \mathbf Z v$ with $(v^{2}) = 0$, 
and $NS(M) + T(M)$ is of co-rank one in $H^{2}(M, \mathbf Z)$.

{\bf 4.} Let $g \in {\rm Bir}\, (M)$ be a bimeromorphic automorphism of $M$. We choose a Hironaka's 
resolution $\pi_{1} : Z \longrightarrow M$ of the indeterminacy $I(g)$ of $g$ (of the source $M$) 
and denote by $\pi_{2} : Z \longrightarrow M$ the induced morphism to the target $M$. Recall that $I(g)$ is the minimum 
(necessarily Zariski closed) subset of $M$ s.t. $g$ is holomorphic over $U := U(g) := M \setminus I(g)$  
and that ${\rm codim}\, I(g) \geq 2$. Let $\{E_{j}\}_{j \in J}$ be the set of exceptional prime divisors of $\pi_{1}$.  
Put $E := \cup_{j \in J}E_{j}$ and $I'(g) := \pi_{2}(E)$. These (with reduced structure) are closed analytic subsets of 
$Z$ and (the target) $M$ respectively. Since $\pi_{1} \vert Z \setminus E : Z \setminus E \simeq M \setminus I(g)$ 
is an isomorphism, we have $g(U) = M \setminus I'(g)$. In particular, $g(U)$ is a Zariski open subset of the target $M$.  

Since $\pi_{1}$ is a successive blow-up in smooth centers, one has: 
$$H^{2}(Z, \mathbf K) = \pi_{1}^{*} H^{2}(M, \mathbf K) \oplus \oplus_{j \in J} \mathbf K [E_{j}]\,\, .$$ 
Here $\mathbf K$ is any ring containing $\mathbf Z$ (mostly, $\mathbf Z$, $\mathbf Q$, $\mathbf R$, or $\mathbf C$).

Recall that $K_{M} = 0$. Then, as divisors, $K_{Z} = \sum_{j \in J}a_{j}E_{j}$ for some positive integers 
$a_{j} > 0$ and $(\pi_{2})_{*}K_{Z} = K_{M} = 0$. Since $M$ is K\"ahler, this implies that the analytic subset 
$I'(g) = M \setminus g(U)$ is of codimension $\geq 2$ as well.

Let $\iota : U \longrightarrow M$ be the natural inclusion. Recall that ${\rm codim}\, M\setminus U \geq 2$ 
and $\Omega_{M}^{2}$ is locally free. 
Then one has $\iota_{*}H^{2}(U, \Omega_{U}^{2}) = 
H^{2}(M, \Omega_{M}^{2}) = \mathbf C \sigma_{M}$. Thus, 
$$g^{*}\sigma_{M} := (\iota_{U})_{*}(g \vert g(U))^{*}(\sigma_{M}\vert g(U))$$  
is a holomorphic $2$-form on $M$. Hence $g^{*}\sigma_{M}$ is a non-zero multiple of $\sigma_{M}$, 
and is in particular everywhere non-degenerate. So, the morphism $g\vert U : U \longrightarrow g(U)$ is an isomorphism 
and we have then $I(g^{-1}) = I'(g)$ as well (cf. [GHJ, Proposition 21.6]). 
Since $\pi_{i}^{*}H^{0}(M, \Omega_{M}^{2}) = H^{0}(Z, \Omega_{Z}^{2})$, we have also 
$\pi_{1}^{*}(g^{*}\sigma_{M}) = \pi_{2}^{*}\sigma_{M}$ as forms. Thus $g^{*}\sigma_{M} = 
(\pi_{1})_{*}\pi_{2}^{*}\sigma_{M}$ 
in $H^{2}(M, \mathbf C)$. In particular, the deinition of $g^{*}$ here is 
compatible with $g^{*}$ on $H^{2}(M, \mathbf Z)$ defined below.  

{\bf 5.} Let us study representations of ${\rm Bir}\, (M)$.
\begin{lemma} \label{lemma:unit}
Let $g \in {\rm Bir}\, (M)$ and define $\chi(g) \in \mathbf C$ by $g^{*}\sigma_{M} = \chi(g)\sigma_{M}$. Then 
$\chi(g)$ is on the unit circle $S^{1}$, i.e. $\vert \chi(g) \vert = 1$. 
\end{lemma}
\begin{proof} Set $\sigma' = g^{*}\sigma_{M}$. Then, by $\pi_{1}^{*}\sigma' = \pi_{2}^{*}\sigma_{M}$,  
$$\int_{M}(\sigma'\overline{\sigma'})^{m} = \int_{Z}(\pi_{1}^{*}\sigma'\overline{\pi_{1}^{*}\sigma'})^{m} 
= \int_{Z}(\pi_{2}^{*}\sigma_{M}\overline{\pi_{2}^{*}\sigma_{M}})^{m} = 
\int_{M}(\sigma_{M}\overline{\sigma_{M}})^{m} = 1\,\, .$$ 
Since $\sigma' = \chi(g)\sigma_{M}$, this impies the result. 
\end{proof}  

The correspondence $Z$ also induces 
a group automorphism 
$g^{*} = (\pi_{1})_{*}\pi_{2}^{*}$ of $H^{2}(M, \mathbf Z)$ (and a $\mathbf K$-linear automorphism on 
$H^{2}(M, \mathbf K)$ as well). By the obsrvation above, $g^{*}$ preserves the 
Hodge decomposition. 
Note that if $g^{*}\eta = \eta'$, then $\pi_{1}^{*}\eta' = \pi_{2}^{*}\eta + \sum_{j \in J}a_{j}E_{j}$ in 
$H^{2}(M, \mathbf K)$ for some $a_{j} \in \mathbf K$. 
It is observed by [GHJ, Section 27] that $(\pi_{1}^{*}\sigma_{M})^{m-1}\vert E_{i} = 0$ for 
each exceptional prime divisor $E_{i}$. Now, by calculating $q_{M}(g^{*}\eta)$ and $q_{M}(\eta)$ 
(by pulling back to $Z$ as above), one obtains:
\begin{theorem} \label{theorem:bir} {\rm [GHJ, Proposition 25.14]}
${\rm Bir}\, (M)$ acts on $H^{2}(M, \mathbf Z)$ as a Hodge isometry with respect 
to BF-form. Moreover, ${\rm Bir}\, (M)$ preserves the positive cone $\mathcal C$. 
\end{theorem}
\begin{proof} Since the last claim is implicit there, we shall give its proof here. Let $\eta$ be a K\"ahler class 
and put $\eta' = g^{*}\eta$. It suffices to show that $(\eta', \eta) \geq 0$. We can write 
$\pi_{1}^{*}\eta' = \pi_{2}^{*}\eta + \sum_{j \in J}a_{j}E_{j}$ with $a_{j} \in \mathbf R$. By using the definition of 
BF-form 
and the fact that $(\pi_{1}^{*}\sigma_{M})^{m-1}\vert E_{i} = 0$, one calculate  
$$(\eta', \eta) = c\int_{M}(g^{*}\eta)\eta(\sigma\overline{\sigma})^{m-1} 
= c\int_{Z}(\pi_{1}^{*}\eta')(\pi_{1}^{*}\eta)\pi_{1}^{*}(\sigma\overline{\sigma})^{m-1}$$ 
$$= c\int_{Z}(\pi_{2}^{*}\eta)(\pi_{1}^{*}\eta)\pi_{1}^{*}(\sigma\overline{\sigma})^{m-1} \geq 0\,\, .$$ 
Here $c$ is a positive constant. The last inequality is because each of three terms (in the last integral) is a pull 
back of a weakly-positive form by a morphism. 
\end{proof} 
 
{\bf 6.} We have then four natural representations of ${\rm Bir}\,(M)$:

$$r : {\rm Bir}\, (M)\, \longrightarrow\, {\rm O}(H^{2}(M, \mathbf Z))\,\, 
;\,\, f \mapsto f^{*}\,\, ,$$

$$\chi : {\rm Bir}\,(M) \longrightarrow 
\mathbf C^{\times}\,\, \text{by}\,\, f^{*}\sigma_{M} = 
\chi(f)\sigma_{M}\,\,  ,$$ 

$$r_{T} : {\rm Bir}\, (M)\, \longrightarrow\, {\rm O}(T(M))\,\, ;\,\, 
f \mapsto f^{*}\vert T(M)\,\, ,$$

$$r_{NS} : {\rm Bir}\, (M) \longrightarrow {\rm O}(NS(M))\,\, ;\,\, 
f \mapsto f^{*}\vert NS(M)\,\, .$$
The next theorem due to Huybrechts [Hu, Section 9] basically reduces our study of ${\rm Bir}\,(M)$ 
to its representation $r({\rm Bir}\,(M))$:
\begin{theorem} \label{theorem:ker} {\rm [Hu]} 
$\vert {\rm Ker}(r : {\rm Bir}\,(M) \longrightarrow {\rm O}(H^{2}(M, \mathbf Z))) \vert < \infty$. 
\end{theorem} 
\begin{proof} This is so important for us that we shall review the proof here. Let $K$ be the kernel. 
Let $g \in K$. Then $g^{*}(\eta) = \eta$ for a K\"ahler class $\eta$. Thus $g \in {\rm Aut}\, (M)$ 
by a result of Fujiki [Fu1] (see also [GHJ, Proposition 27.7] for a useful refinement). 
Recall that by Yau's solution 
of Calabi's conjecture [Ya1,2], there is a unique Ricci-flat $\tilde{\eta}$ 
metric associated to the class $[\eta]$. Thus, we have also $K < {\rm O}(M, \tilde{\eta})$, 
the group of isometries of $(M, \tilde{\eta})$. 
This ${\rm O}(M, \tilde{\eta})$ is a compact subgroup of ${\rm Diff}(M)$ (with respect to the compact-open topology), 
because 
$M$ is compact (see for instance [He] or [Kb]). On the other hand, ${\rm Aut}\,(M)$ forms a closed discrete subgroup 
of ${\rm Diff}(M)$ by $H^{0}(M, T_{M}) = 0$. 
Thus, $\vert {\rm Aut}\,(M) \cap {\rm O}(M, \tilde{\eta}) \vert < \infty$ and hence 
$\vert K \vert < \infty$ as well. 
\end{proof}

The next Theorem slightly strengthens earlier results of 
Nikulin [Ni] and Beauville [Be2, Proposition 7]: 
  
\begin{theorem} \label{theorem:conv} 
Let $M$ be a hyperk\"ahler manifold. Let $G < {\rm Bir}\, (M)$. 
Set $A' := \chi (G)$ and $A := r_{T}(G)$. Then: 
\begin{list}{}{
\setlength{\leftmargin}{10pt}
\setlength{\labelwidth}{6pt}
}
\item[(1)] The characteristic polynomial $\Phi(x)$ of $g^{*} \vert T(M)$ ($g \in G$)  
is primary in $\mathbf Z[x]$, i.e. $\Phi(x)$ is of the form $f(x)^{m}$ for an irreducible polynomial 
$f(x) \in \mathbf Z[x]$. 
\item[(2)] The natural homomorphism 
$$\psi : A \longrightarrow A'\,\, ;\,\, g^{*} \vert T(M) \mapsto \chi(g)$$ 
is an isomorphism. In particular, $A$ is an abelian group. 
\item[(3)] $A' < S^{1}$. 
\item[(4)] If $NS(M)$ is hyperbolic, then $A \simeq A'$ is 
a finite cyclic group $\mu_{n}$, where $n \geq 1$ is some integer s.t. $\varphi(n) \vert {\rm rank}\, T(M)$. 
\item[(5)] If $NS(M)$ is parabolic, then $A \simeq A' = \{1\}$. 
\item[(6)] If $NS(M)$ is elliptic, then $A \simeq A'$ is torsion free, i.e. 
either $A' = \{1\}$ or ${\rm ord}\,(\chi(g)) = \infty$ for all $g \in G \setminus \{1\}$. Moreover, 
${\rm Ker}\,(\chi : G \longrightarrow A')$ is finite, and $\vert G \vert = \infty$ iff $\vert A' \vert = \infty$. 
\end{list}
\end{theorem} 

\begin{proof}
Let us show (1). If otherwise, we have a decomposition in $\mathbf Z[x]$:
$$\Phi(x) = f_{1}(x)f_{2}(x)\,\, {\rm s.t.}\,\, (f_{1}(x), f_{2}(x)) = 1\,\, .$$ 
This gives a decomposition of $T(M)_{\mathbf Q}$, say $T(M)_{\mathbf Q} = V_{1} \oplus V_{2}$ 
such that both $V_{i}$ are $g^{*}$-stable and the characteristic polynomial of $g^{*} \vert V_{i}$ 
is $f_{i}(x)$ for each $i = 1$, $2$. Since $f_{1}$ and $f_{2}$ have no common zero, $\sigma_{M}$ 
then belongs to one of $V_{i} \otimes 
\mathbf C$, a contradiction to the minimality of $T(M)$. 

Let us show (2). Clearly $\psi$ is surjective. 
Take $g \in G$ such that 
$\chi(g) = 1$, i.e. $g^{*}\sigma_{M} = \sigma_{M}$. Let 
$V := T(M)^{g^{*}}$ be the set of $g$-invariant elements of 
$T(M)$. Then $\sigma_{M} \in V \otimes \mathbf C$, whence, by the minimality of $T(M)$, 
we have $V = T(M)$, i.e. $g^{*} \vert T(M) = id$. The last statement of (2) 
is now clear, because so is $A' < \mathbf C^{\times}$. The assertion (3) is proved in Lemma (2.1).

Let us show (4). We first show that $\chi(g)$ is 
a root of unity for $g \in G$. Set 
$$N := P_{T(M) \otimes \mathbf R}^{\perp} := 
\{v \in T(M) \otimes \mathbf R\, \vert\, (v, P) = 0\}\,\, ,$$ 
where $P = \mathbf R \langle {\rm Re}\,\sigma_{M}\, ,\, {\rm Im}\, \sigma_{M} \rangle$ is the positive $2$-plane. 
Then $T(M) \otimes \mathbf R = P \oplus N$ and $N$ is negative definite. Put $g_{T} := g^{*} \vert T(M)$. 
Then $g_{T}(P) = P$ and $g_{T}(N) = N$. Hence 
$g_{T} \in {\rm O}\,(P) \times {\rm O}\,(N)$. Thus the eigenvalues of 
$g_{T}$ are on the unit circle. On the other 
hand, the eigenvalues of $g_{T}$ 
are all algebraic integers, because $g_{T} \in {\rm O}\,(T(M))$ so that 
$g_{T}$ is represented by an integral matrix. 
Hence, they are roots of unity by 
Kronecker's theorem (see eg. [Ta, Section 9.2]). 
Since $\chi(g)$ is one of the eigenvalue of $g_{T}$, it is also a root of unity, 
say, a primitive $n(g)$-th root of unity. By (1), $\varphi(n(g)) \vert {\rm rank}\, T(M)$. Hence the orders $n(g)$ 
of $g$ ($g \in G$) are bounded. Thus $A' \simeq \mu_{n}$ for some $n$ with $\varphi(n) \vert {\rm rank}\, T(M)$.  

Let us show (5). Set $NS(M) \cap T(M) = \mathbf Z v$ (as before). Replacing $v$ by $-v$ if necessary, 
we may assume that $v \in \partial \mathcal C$. 
Then, by $g^{*}(\mathcal C) = \mathcal C$, we have $g^{*}(v) = v$ for each $g \in G$. Thus 
$g^{*}\sigma_{M} = \sigma_{M}$ for each $g \in G$ by (1), i.e. $A' = \{1\}$. 

Finally, we shall show (6). Since $NS(M) \oplus T(M)$ is of finite index in $H^{2}(M, \mathbf Z)$, the group 
${\rm Ker}\, r_{T}/{\rm Ker}\, r$ is embedded into ${\rm O}\,(NS(M))$. 
Thus ${\rm Ker}\, \chi = {\rm Ker}\, r_{T}$ is finite by $\vert {\rm O}\,(NS(M)) \vert < \infty$ 
and by Theorem (2.2). Thus $\vert G \vert < \infty$ iff $\vert A' \vert < \infty$. 
It remain to show that $\chi(g) = 1$ if ${\rm ord}\, \chi(g) < \infty$. By the previous argument, the group $H := \langle g \rangle$ is finite.   
Note that $H^{2}(M, \mathbf Q)$ is dense in $H^{2}(M, \mathbf R)$ 
and the positive 
cone $\mathcal C$ is open and non-empty in $H^{1,1}(M, \mathbf R)$. Then, since 
the 
natural projection 
$p : H^{2}(M, \mathbf R) \longrightarrow H^{1,1}(M, \mathbf R)$ is continuous, 
there is $\eta \in H^{2}(M, \mathbf Q)$ such that 
$\eta^{(1,1)} := p(\eta) \in \mathcal C$. Since $\mathcal C$ 
is a convex cone and stable under 
${\rm Bir}\, (M)$, it follows that the element 
$$\tilde{\eta} := \sum_{h \in H} h^{*}\eta \,\, $$ 
is $H$-invariant, rational and satisfies 
$\tilde{\eta}^{(1,1)} \in \mathcal C$. Assume that $\chi(g) \not= 1$. Then, 
$\tilde{\eta}$ would have no $(2,0)$ and $(0, 2)$ 
component, because it is $H$-invariant. That is, $\tilde{\eta}$ would be  
of pure $(1,1)$ type. Since $\tilde{\eta}$ is rational, we would then have 
$\tilde{\eta} \in NS(M) \otimes \mathbf Q$ by the Lefschetz (1,1)-Theorem. 
But $(\tilde{\eta}^{2}) > 0$ by the choice of $\tilde{\eta}$, 
a contradiction to the negative definiteness of $NS(M)$. Thus $\chi(g) = 1$. 

\end{proof}

This theorem (together with the following examples) shows that the value $\chi(g)$ is an effective invariant 
which can distinguish 
certain non-projective hyperk\"ahler manifolds, especially those with elliptic N\'eron-Severi lattices, from projective 
ones. This view point will be important in Section 8. 
\begin{example} \label{example:pair} 
Let $S$ be a K3 surface having an automorphism $f$ with a Siegel disk $D$, 
say, $(D, f\vert D)$ is isomorphic to 
$(\Delta^{2}, F)$ with $F^{*}(z_{1}, z_{2}) = 
(\alpha_{1} z_{1}, \alpha_{2} z_{2})$. Here $\alpha_{1}$ 
and $\alpha_{2}$ are multiplicatively independent. Then,   
$\chi(f) = \alpha_{1} \alpha_{2}$ and $\chi(f)$ is not a root of unity. 
So, $NS(S)$ is elliptic by Theorem (2.3) and therefore the 
algebraic dimension $a(S) = 0$ by the classification theory of surfaces (cf. [BPV]).
\end{example} 

We shall give examples of hyperk\"ahler manifolds with similar automorphisms:
 
\begin{example} \label{example:beauville}
Let $m \geq 2$ be an integer and $(S, f)$ be a McMullen's pair. 
Let $(M_{m}, f_{m})$ be the pair of the Duady space $M_{m} := S^{[m]}$ 
(cf. Introduction) and its automorphism $f_{m}$ naturally induced 
by $f$. Then, $M_{m}$ is a $2m$-dimensional hyperk\"ahler manifold 
such that $NS(M_{m}) \simeq \langle -2(m-1) \rangle$. In particular, $NS(M_{m})$ 
is elliptic and $\rho(M_{m}) = 1$. Moreover, $\chi(f_{m})$ is not 
a root of unity. 
\end{example} 

\begin{proof} The first part is shown by [Be1]. The fact $NS(M_{m}) 
\simeq \langle -2(m-1) \rangle$ follows from $NS(S) = \{0\}$ and [ibid]. 
By [ibid], a symplectic from $\sigma_{M_{m}}$ of $M_{m}$ is 
given as follows: the $2$-from $\sum_{i=1}^{m} \text{pr}^{*}_{i}\sigma_{S}$ 
on $\Pi_{1}^{m} S$ 
descends to the $2$-form $\sigma_{S^{(m)}}$ on the Chow variety
$S^{(m)} := (\Pi_{i=1}^{m} S)/S_{m}$. The pullback of $\sigma_{S^{(m)}}$ 
under the Hilbert-Chow morphism $M_{m} \longrightarrow S^{(m)}$ 
gives $\sigma_{M_{m}}$.

By definition of 
$M_{m}$, the automorphism $f$ of $S$ naturally 
induces an automorphism $f_{m}$ of $M_{m}$. By the description of 
$\sigma_{M_{m}}$, 
one has   
$f_{m}^{*}\sigma_{M_{m}} = \chi(f)\sigma_{M_{m}}$.  
Thus $\chi(f_{m}) = \chi(f)$ and is not a root of unity by Example (2.4). 
\end{proof} 

The following corollary says that geometric action 
of $G (< {\rm Bir}\, (M))$ on $M$ can be approximated by its representation on $NS(M)$ 
when $NS(M)$ is hyperbolic or parabolic.

\begin{corollary} \label{corollary:nen} 
Let $M$ be a hyperk\"ahler manifold. Let $G < {\rm Bir}\, (M)$. 
Let $K := {\rm Ker}\, (r_{NS} : G \longrightarrow {\rm O}\, (NS(M)))$ and $H := r_{NS}(G)$. 
Assume that $NS(M)$ is hyperbolic or parabolic. Then $K$ 
is finite. In particular, $G$ is almost abelian 
iff so is $H$, and they have the same rank. 
\end{corollary} 

\begin{proof} Let $G^{*} := r(G) = {\rm Im}\,(r : G \longrightarrow {\rm O}(H^{2}(M, \mathbf Z)))$. 
By Theorem (2.2) and Proposition (9.3) in the appendix, it suffices to show that the group 
$K^{*}\, :=\, \{g \in G^{*}\,\, \vert\,\, g \vert NS(M) = id\}$ is finite.

Consider the natural map $\tau : K^{*}\, \longrightarrow\, r_{T}(G)\,\, ;\,\, g\, \mapsto\, g 
\vert T(M)$.

By Theorem (2.3)(4)(5), $\vert r_{T}(G) \vert < \infty$. So, it suffices to 
show that $\tau$ is injective.

First consider the case where $NS(M)$ is hyperbolic. 
Then, $NS(M) \oplus T(M)$ is 
of finite index in $H^{2}(M, \mathbf Z)$. Thus, $\tau$ 
is injective.

Next consider the case where $NS(M)$ is parabolic. 
The remaing argument is quite similar to 
that of [Og1, Appendix].

By the assumption, we have $NS(M) \cap T(M) = \mathbf Z v$. This $v$ is 
primitive 
and satisfies $(v^{2}) = 0$ and $v^{\perp}_{NS(M)} = NS(M)$ 
and $v^{\perp}_{T(M)} = T(M)$. 
Then, BF-form on $NS(M)$ descends to the bilinear 
form on  $\overline{NS}(M) := NS(M)/\mathbf Z v$ and makes $\overline{NS}(M)$ 
a negative definite lattice of rank $\rho(M)-1$. 
Similarly, BF-form on $T(M)$ descends to the bilinear 
form on $\overline{T}(M) := T(M)/\mathbf Z v$ and makes $\overline{T}(M)$ 
a non-degenerate lattice. The isometry $G^{*}$ also descends to 
the isometry of both $\overline{NS}(M)$ 
and $\overline{T}(M)$.

Let $\langle v\, ,\, u_{i}\, (1 \leq i \leq \rho(M) -1)\,\rangle$ 
be an integral 
basis of $NS(M)$ 
and $\langle v\, ,\, t_{j}\, (1 \leq j \leq k)\, \rangle$ be an integral 
basis of $T(M)$. Let $w \in H^{2}(M, \mathbf Z)$ be an element such that 
$$\langle v\, ,\,  t_{j}\, (1 \leq j \leq k)\, ,\,  u_{i}\, (1 \leq i \leq 
\rho(M)-1)\, , w \rangle$$
forms a basis of $H^{2}(M, \mathbf Q)$. Let $g \in {\rm Ker}\, \iota$. Then we have 
$$g(v) = v\,\, ,\,\, g(u_{i}) = u_{i}\,\, , \,\, 
g(t_{j}) = t_{j}\,\,  ,$$ 
$$g(w) = a(g)w + b(g)v + \sum_{i=1}^{\rho(M)-1} c_{i}(g)u_{i} + 
\sum_{j=1}^{k} d_{j}(g)t_{j}\,\,  .$$
Here $a(g)$, $b(g)$, $c_{i}(g)$ and $d_{j}(g)$ are rational numbers. 

Note that $(v^{2}) = (v, u_{i}) = (v, t_{j}) = 0$. Hence $(v, w) \not= 0$, 
because $H^{2}(X, \mathbf Q)$ is non-degenerate. 
Thus $a(g) = 1$ from $(w, v) = (g(w), g(v))$. From $(g(w), g(u_{l})) = 
(w, u_{l})$, 
we obtain 
$$\sum_{i=1}^{\rho(M)-1}c_{i}(g)(u_{i}, u_{l}) \,\, = \,\, 0\,\,$$ 
for all $l$ with $1 \leq l \leq \rho(M) -1$. Note that $\langle u_{i}\, 
\text{mod}\, \mathbf Z v \rangle_{i=1}^{\rho(M)-1}$ forms 
integral basis of $\overline{NS}(M)$ and recall 
that $\overline{NS}(M)$ is of negative 
definite. Thus, we have 
$c_{i}(g) = 0$ for all $i$ from the equality above. Using the fact that 
$\langle t_{j}\, \text{mod}\, \mathbf Z v \rangle_{j=1}^{k}$ forms 
integral basis 
of $\overline{T}(M)$ and $\overline{T}(M)$ is non-degenerate, 
we also obtain that $d_{j}(g) = 0$ for all $j$ in a similar manner. 
Thus $g(w) = w + b(g)v$. Since $(g(w)^{2}) = (w^{2})$ and $(w, v) \not= 0$, 
we have finally $b(g) = 0$. Hence $g(w) = w$. Thus 
$g = id_{H^{2}(M, \mathbf Q)}$. Therefore $\tau$ is injective.  
\end{proof}

\section{Salem polynomials and transcendental representations} 

Throughout this section, $M$ is a hyperk\"ahler manifold with elliptic $NS(M)$. 
We will study the transcendental representation of 
${\rm Bir}\,(M)$ more closely by means of the notion of Salem polynomials.
Our main result of this Section is Theorem (3.4).

By Kronecker's theorem (see eg. [Ta, Section 9.2]), a monic 
irreducible polynomial in $\mathbf Z[t]$ is {\it cyclotomic} iff all the roots 
are on the unit circle $S^{1}$. Salem polynomials form the second simplest 
class of 
monic irreducible polynomials in $\mathbf Z[t]$ in the view: {\it "how 
many roots are off the unit circle ?"}. Precise definition is:   

\begin{definition} \label{definition:salem} (cf. [Mc], [BM])
An {\it irreducible monic} polynomial $\Phi(t) \in \mathbf Z[t]$ of 
degree $n$ is a {\it Salem 
polynomial} if the roots of $\Phi(t)$ consist of 
two real roots $\alpha$ and $1/\alpha$ such that $\alpha > 1$ and 
$n-2$ roots on the unit circle. The unique root $\alpha >1$ is 
the {\it Salem number} associated with 
$\Phi(t)$.   
\end{definition}  

\begin{remark} \label{remark:salemdeg} $\text{deg}\, \Phi(t) 
\equiv 0\, (\text{mod}\, 2)$ and $\Phi(0) = 1$ if $\Phi(t)$ is a 
Salem polynomial. 
\end{remark}
Recall that cyclotomic polynomials with fixed degree 
are at most finitely many (possibly empty). On the other hand, as it is 
shown by Gross and McMullen [GM], there are infinitely many Salem 
polynomials in each even degree $\geq 6$. The following finiteness criterion about 
Salem polynomials is in [Mc, Section 10]:

\begin{proposition} \label{proposition:finsalem} 
There are at most finitely many Salem polynomials with fixed degree and 
bounded trace. That is, given an integer $n > 0$ 
and real numbers 
$B_{1} < B_{2}$, there are at most finitely many Salem polynomials 
$\varphi(t)$ such that 
$${\rm deg}\, \varphi(t) = n\,\, {\rm and}\,\, B_{1} < {\rm tr}\, 
\varphi(t) < B_{2}\,\, .$$ 
Here the trace of $\varphi(t)$ 
is the sum of roots counted with multiplicities, and is the same as the minus 
of the coefficient of the second leading term $t^{n-1}$. 
\end{proposition}

\begin{proof} The argument here is taken from [Mc, Section 10]. 
Write 
$$\varphi(t) = t^{n} + \sum_{k=1}^{n-1}a_{k}t^{k} + 1\,\, .$$
Let $\alpha > 1$, $1/\alpha$, $\beta_{i}$, 
$\overline{\beta_{i}}$ 
($1 \leq i \leq (n-1)/2$) be the roots of $\varphi(t)$. Since 
$\beta_{i}$ are on the unit circle, one has

$$-(n-1) \leq \sum_{i=1}^{(n-1)/2} \beta_{i} + 
\overline{\beta_{i}} \leq n-1\,\, .$$ 
On the other hand, by the assumption, one has  
$$B_{1} < \text{tr}\, \varphi(t) = \alpha + 1/\alpha +  
\sum_{i=1}^{(n-1)/2} \beta_{i} 
+ \overline{\beta_{i}} < B_{2}\,\, .$$ 
Thus, the Salem number $\alpha > 1$ is also bounded from the above. 
Hence the roots of $\varphi(t)$ are in some compact set of $\mathbf C$ 
(independent of $\varphi(t)$). Up to sign, $a_{k}$ is the elementary 
symmetric function of degree $k$ of the roots. 
Thus $a_{k}$ are also in some compact set 
of $\mathbf C$ (independent of $\varphi(t)$ and $k$). 
Since $a_{k}$ are all integers, this implies the result. 
\end{proof}
 
\par \vskip 1pc

As it is explained in Introduction, McMullen's pairs are closely related 
to Salem polynomials. Salem polynomials also play an important role in our 
study:

\begin{theorem} \label{theorem:nonproj} 
Let $M$ be a hyperk\"ahler manifold with elliptic $NS(M)$. 
Let $f \in {\rm Bir}\,(M)$. Set $F := f^{*} \vert T(M)$. 
If ${\rm ord}\, F = \infty$ (or equivalently ${\rm ord}\, f = \infty$ by 
Theorem (2.4)(6)), then the characteristic polynomial 
$$\Phi(t) := \text{det}\, (tI - F)$$ 
of $F$ is a Salem polynomial. In particular, 
${\rm rank}\, T(M) \equiv 0\, {\rm mod}\, 2$ and 
$f^{*} \vert T(M) \in {\rm SO}\,(T(M))$, if there is $f \in {\rm Bir}\, (M)$ 
with ${\rm ord}\, f = \infty$. 
\end{theorem} 

\begin{proof} Recall that a complex number $\lambda$ is 
an eigenvalue of 
$F$ iff $\Phi(\lambda) = 0$. We denote by $E(\lambda)$ the generalized 
$\lambda$-eigenspace of $F$, that is, $\Phi(\lambda) = 0$ and  
$$E(\lambda) := 
\{v \in T(M) \otimes \mathbf C\, \vert\, (\lambda I - F)^{n}v = 0\, 
\, \exists n \in \mathbf Z_{>0} 
\}\,\, .$$ 
We proceed our proof by dividing into several steps. 

\begin{claim} \label{claim:nonunity} $\chi(f)$ is not a root of unity.  
\end{claim} 

\begin{proof} Otherwise, $\chi(f) = 1$ by Theorem (2.4)(6).  
Then $F = id$ by Theorem (2.4)(2), a contradiction. 
\end{proof} 

\begin{claim} \label{claim:big} $\Phi(t)$ has at least one 
root $\alpha_{0}$ {\it outside} the unit circle $S^{1}$, that is, 
there is a complex number $\alpha_{0}$ such that $\Phi(\alpha_{0}) = 0$ 
and $\vert \alpha_{0} \vert > 1$.  
\end{claim} 

\begin{proof} Note that $\text{det}\, F = \pm 1$ by 
$F \in \text{O}(T(M))$ and by the fact that $T(M)$ is non-degenerate 
if $NS(M)$ is elliptic. So, otherwise, the roots of $\Phi(t)$ would be all on 
$S^{1}$. On the other hand, they 
are algebraic integers. Hence all the roots of $\Phi(t)$ would the 
be roots of unity by Koronecker's theorem. In particular, 
$\chi(f)$ would be a root of unity, a contradiction to Claim (3.5). 
\end{proof} 

The next Lemma is a special case, the case $p =1$, of a well-known 
fact about the real orthogonal group $\text{O}\,(p, q)$ (cf. [Mc, Lemma 3.1]): 

\begin{lemma} \label{lemma:hodgeindex} Let $g \in {\rm O}\,(1, m)$. Then 
$g$ has at most one 
eigenvalue (counted with multiplicity) outside $S^{1}$. 
\end{lemma} 

\begin{claim} \label{claim:big} $\Phi(t)$ has exactly one real root 
$\alpha$ such that $\alpha > 1$, and exactly two real roots, 
which are $\alpha$ and $1/\alpha$, counted with multiplicities. 
The other roots are on the unit circle. 
\end{claim} 

\begin{proof} 
Set $S := \oplus_{\vert \lambda \vert > 1} E(\lambda)$. There is an 
$F$-stable $\mathbf R$-linear 
subspace $R \subset T(M) \otimes 
\mathbf R$ such that $S = R \otimes \mathbf C$; Indeed, if 
$\lambda$ is an eigenvalue 
of $F$ with $\vert \lambda \vert > 1$ and $v \in E(\lambda)$, then so is 
$\overline{\lambda}$ and $\overline{v} \in E(\overline{\lambda})$.

Choose $\alpha_{0}$ as in Claim (3.6). One has 
$0 \not= E(\alpha_{0}) \subset S$. Thus $R \not= \{0\}$. 
We show $\text{dim}_{\mathbf R} R \leq 1$. Let $P$ be the positive $2$-plane 
in $T(M) \otimes \mathbf R$ and $N$ the orthogonal complement of $P$ in $T(M) \otimes \mathbf R$. 
Then $R \subset N$ by $F(P) = P$, $F(N) = N$ and by the fact that 
$\chi(f)$ and $\overline{\chi(f)}$, the eigenvalues 
of $F \vert P$, are on $S^{1}$. On the other hand, since 
the signature of $N$ is $(1,0, t)$, one has $\text{dim}_{\mathbf R}\, R \leq 1$ by Lemma (3.7). 
So, $\text{dim}_{\mathbf R}\, R = 1$ by $R \not= \{0\}$. Thus $\Phi(t)$ has 
exactly one root $\alpha$ such that $\vert \alpha \vert > 1$ and that 
this $\alpha$ is necessarily real and of multiplicity one.

Next we show that $\alpha > 0$. Let $\mathcal C$ be the positive cone of $M$. 
Since $\alpha \in \mathbf R$, there is 
$0 \not= v \in N$ such that $F(v) = \alpha v$. Note that 
$v \in H^{1,1}(M, \mathbf R)$ by $N \subset H^{1,1}(M, \mathbf R)$. 
One calculates $(v, v) = (F(v), F(v)) = \alpha^{2}(v, v)$. 
Then $(v, v) = 0$ by $\vert \alpha \vert > 1$. So, replacing $v$ by $-v$ 
if necessary, one has $v \in \partial\, \mathcal C$. Here 
$\partial\, \mathcal C$ is the boundary of the positive cone. 
Since $f^{*}(\mathcal C) 
= \mathcal C$, we have then $\alpha > 0$. Thus $\alpha > 1$ 
by $\vert \alpha \vert > 1$.

Since the same holds for $f^{-1}$, it follows that $\Phi(t)$ has exactly one root 
$\beta$ {\it inside} $S^{1}$, 
and this $\beta$ is real, of multiplicity one and satisfies $0 < \beta < 1$. 
In particular, the remaining roots ($\not= \alpha, \beta$) of $\Phi(t)$ 
are on $S^{1}$. Now $\beta = 1/\alpha$ 
by $\text{det}\, F = \pm 1$, and we are done. 
\end{proof}  

\begin{claim} \label{claim:irred} $\Phi(t)$ is 
irreducible in $\mathbf Z [t]$.
\end{claim} 

\begin{proof} Since $\Phi(t)$ is monic, it suffices to show that $\Phi(t)$ 
is irreducible in $\mathbf Q[t]$. By Claim (3.8), there is no $\varphi(t) \in \mathbf Q[t]$ such that 
$\Phi(t) = (\varphi(t))^{m}$ for $\exists m \geq 2$. So, the result follows 
from Theorem (2.4)(1). 
\end{proof} 

The next claim completes the proof of Theorem (3.4).  

\begin{claim} \label{claim:orth} $F \in \text{SO}(T(M))$.  
\end{claim} 

\begin{proof} It is clear that $F \in \text{O}(T(M))$. It suffices 
to show $\text{det}\, F = 1$. We already know that 
$\Phi(t) \in \mathbf Z[t]$ is an irreducible polynomial 
of $\text{deg}\, \Phi(t) = \text{rank}\, T(M) \geq 2$. Thus $\Phi(\pm 1) \not= 0$. So $\Phi(t)$ has no real 
root on the unit circle. 
Thus, the roots other than $\alpha$, $1/\alpha$ are all imaginary, 
and appear in pairs like $\lambda$, $\overline{\lambda}$. Hence 
the product of all the roots of $\Phi(t)$ is $1$ and 
$\text{det}\, F = 1$.  
\end{proof}
Now we are done. Q.E.D. of Theorem (3.4).  
\end{proof}

\begin{corollary} \label{corollary:simple} Under the same assumption as in 
Theorem (3.4), each eigenvalue of $F$ is 
of multiplicity one, and $F$ is diagonalizable on $T(M) \otimes \mathbf C$.   
\end{corollary}

\begin{proof} This is because $\Phi(t)$ is irreducible. \end{proof}

\begin{remark} \label{remark:fingen} By Theorem (2.4)(6), the main part of Theorem (1.5)(1) follows 
once we will know ${\rm Im}\, \chi \simeq \mathbf Z$ when it is infinite. However, at the moment, 
even the finite generation of 
${\rm Im}\, \chi$ is not so clear. For instance, 
$S^{1} (< \mathbf C^{\times})$ has a subgroup isomorphic to the additive group 
$\mathbf Q$, 
namely 
$$\{w \in \mathbf C^{\times}\, \vert\, w =\text{exp}\, 
(2\pi \sqrt{-1}\sqrt{2}q), q 
\in \mathbf Q\}\,\, .$$ 
The additive group $\mathbf Q$ is not finitely generated and 
any two elements of $\mathbf Q$ are not linearly 
independent over $\mathbf Z$. Also, 
$S^{1}$ has a subgroup isomorphic to the additive group 
$\mathbf Z^{r}$ for each $r > 0$. For instance, one 
can check that the following subgroup is isomorphic to $\mathbf Z^{r}$:  
$$\{w \in \mathbf C^{\times}\, \vert\, w =\text{exp}\, 
(2\pi\sqrt{-1} (\sum_{k=1}^{r}a_{k}n_{k}))\,\, ,\,\, (n_{k})_{k=1}^{r}\, 
\in\, \mathbf Z^{r} \}\,\, ,$$
where $a_{k} := \text{log}_{2}\, p_{k}$ and $p_{k}$ are mutually different 
odd prime numbers. We will study finite generation of ${\rm Im}\, \chi$ 
and its rank in the next section. 
\end{remark}

\section{automorphism group of a hyperk\"ahler manifold 
with elliptic N\'eron-Severi lattice}

In this section, we shall prove the main part of Theorem (1.5)(1) deviding into several steps. 
Existence part in dimension $2$ will be shown in Section 6. Cruicial steps are Propositions (4.3) and (4.4).

Let $M$ be a hyperk\"ahler manifold with elliptic $NS(M)$. 
Put: 
$$A' := \text{Im}\, (\chi : {\rm Bir}\,(M) \longrightarrow 
\mathbf C^{\times})\,\, \text{and}\,\, 
A := \text{Im}\, ({\rm Bir}\,(M) \longrightarrow SO(T(M)))\,\, .$$ 
By Theorem (2.4)(2), there is a natural isomorphism: 
$$\psi : A \simeq A'\,\, ;\,\, g^{*} \vert T(M) 
\mapsto \chi(g)\,\, .$$  
As we remarked in (3.12), we may show that $A' \simeq \mathbf Z$ assuming that there is
$f \in {\rm Bir}\,(M)$ such that 
$$\text{ord}\, F = \infty\,\,  ,\,\, \text{where}\,\,  F := f^{*} \vert T(M)\,\, .$$ 
Let $\Phi(t)$ be the characteristic 
polynomial of $F$. We already know that $\Phi(t)$ is a Salem 
polynomial (Theorem (3.4)).  
Thus $\text{deg}\, \Phi(t) = \text{rank}\, T(M)$ must be even. This already 
proves the second statement of (1.5)(1).

Let us compare $F$ 
with another element $G := g^{*} \vert T(M) \in A$ given by 
$g \in {\rm Bir}\,(M)$. 

\begin{claim} \label{claim:comm} 
There is a polynomial $\varphi(t) \in \mathbf Q[t]$ such that 
$G = \varphi(F)$. 
\end{claim} 

\begin{proof} Choose an integral basis of $T(M)$. 
Then, we may regard $F, G \in M(n, \mathbf Q)$, where 
$n := \text{rank}\, T(M)$. 
Since $A$ is an abelian group, one has $FG = GF$. 
Thus $G$ is an element of the $\mathbf Q$-linear subspace $V$
of $M(n, \mathbf Q)$. Here 
$$V := \{X \in M(n, \mathbf Q) \vert FX = XF\,\}\,\, .$$ 
The equation $FX = XF$ is a system of linear homogeneous equations with 
rational coefficients. One can solve this by Gauss' elimination method, 
which obviously commutes with field extensions. Thus   
$$V \otimes \mathbf C = \{X \in M(n, \mathbf C) \vert FX = XF\}\,\, .$$
By Corollary (3.11), there is a basis $\langle v_{k} \rangle_{k=1}^{n}$ 
of $T(M) \otimes \mathbf C$ under which $F$ is represented by a diagonal 
matrix $F = \text{diag}(a_{1}, \cdots , a_{n})$ with 
$a_{i} \not= a_{j}$ for $i \not= j$. Identify $X$ with its matrix 
representation with respect to $\langle v_{k} \rangle_{k=1}^{n}$. 
Now by an explicit matrix calculation, one sees $V \otimes \mathbf C = D(n)$.  
Here $D(n)$ is the set of diagonal matrices of size $n$. 
Thus $\text{dim}_{\mathbf C} V \otimes \mathbf C = n$, 
and hence $\text{dim}_{\mathbf Q} V  = n$. 
Since $\Phi(t)$ is irreducible in $\mathbf Q[t]$, the $n$-elements 
$A^{k}$ ($0 \leq k \leq n-1$) of $V$ are linearly independent over $\mathbf Q$ 
and then form a $\mathbf Q$-basis of $V$. Since $G \in V$, one has then 
$G = \varphi(F)$ for $\exists \varphi(t) \in \mathbf Q[t]$ (with $\text{deg}\, \varphi(t) \leq n-1$). 
\end{proof} 

\begin{claim} \label{claim:unitelem} Let $K (\subset \mathbf C)$ 
be the minimal splitting field of $\Phi(t)$ and $O_{K}$ be the ring of 
integers of $K$, i.e. 
the normalization of $\mathbf Z$ in $K$, and $U_{K} \subset O_{K}$ be 
the group of units of 
$O_{K}$. Then, $\chi(g) \in U_{K}$ for $\forall g \in {\rm Aut}\,(M)$.  
\end{claim} 

\begin{proof} $\chi(g)$ and $\overline{\chi(g)}$ are eigenvalues of 
$G$. Thus, by Claim (4.1), 
$\chi(g)$, $\overline{\chi(g)} \in K$. $\chi(g)$ 
and $\overline{\chi(g)}$ 
are also roots of the monic polynomial $\text{det}(tI - G) 
\in \mathbf Z[t]$. Thus $\chi(g)$, $\overline{\chi(g)} \in O_{K}$. 
By Theorem (2.4)(3), we have also 
$1 = \vert \chi(g) \vert^{2} = \chi(g) \overline{\chi(g)}$. 
Hence $\chi(g)$, $\overline{\chi(g)} \in U_{K}$. 
\end{proof} 

\begin{proposition} \label{proposition:dirichlet} $A$ is a finitely 
generated free 
abelian group.
\end{proposition} 

\begin{proof} By Claim (4.2), one has $A' < \mathbf C^{\times} \cap U_{K} < U_{K}$. 
By Dirichlet's unit theorem (see eg. [Ta, Section 9.3]), $U_{K}$ is 
a finitely generated abelian group. Hence, so is its 
subgroup $A'$ by the fundamental theorem of finitely 
generated abelian groups. Moreover, $A'$ is torsion free by Theorem (2.4)(6). 
Now the result follows from $A \simeq A'$ (Theorem (2.4)(2)). 
\end{proof}

Set $r := \text{rank}\, A$.  The next proposition completes the proof of 
the main part of Theorem (1.5)(1): 

\begin{proposition} \label{proposition:rankr} $r = 1$ (if $r > 0$). 
\end{proposition} 

\begin{proof}  Assuming that $r \geq 2$, we shall derive a contradiction. 
Take a free basis $\langle F_{i} := f_{i} \vert T(M) \rangle_{i=1}^{r}$ of $A$. We may now re-choose $f$ so that $f = f_{1}$. 
Since $r \geq 2$, there is $f_{2}$.  
Set $g := f_{2}$. As before, we write $F = f^{*} \vert T(M)$ and 
$G = g^{*} \vert T(M)$. 
One has $\text{ord}\, F = \text{ord}\, G = \infty$, because $F$ and $G$ form a part of free basis of $A$. Thus, the characteristic polynomials 
$\Phi_{F}(t)$ 
and $\Phi_{G}(t)$ of $F$ and $G$ are both Salem polynomials (Theorem (3.4)). 
Let $\alpha > 1$ and $\beta > 1$ be the Salem numbers of 
$\Phi_{F}(t)$ and $\Phi_{G}(t)$.

Let $(n, m) \in \mathbf Z^{2} - \{(0,0)\}$. Then 
$\text{ord}\, F^{n}G^{m} = \infty$, because $F$ and $G$ form a part of 
free basis of $A$. Thus, the characteristic polynomial 
$\varphi^{(n, m)}(t)$ of 
$F^{n}G^{m}$ is also a Salem polynomial (Theorem (3.4)). In particular, 
$\varphi^{(n, m)}(t)$ has exactly two roots off the unit circle and they are both real. On the other hand, $FG = GF$ and $F$ and $G$ are both 
diagonalizable by Corollary (3.11). Thus $F$ and $G$ are simultaneously 
diagonalizable. Hence the pair of the two (necessarily real) roots of 
$\varphi^{(n,m)}(t)$ off the unit circle must be either 
$$\{\alpha^{n}\beta^{m}\, ,\, \alpha^{-n}\beta^{-m}\}\,\, \text{or}\,\,   
\{\alpha^{n}\beta^{-m}\, ,\, \alpha^{-n}\beta^{m}\}\,\, .$$ 
By replacing $g = f_{2}$ by $g = f_{2}^{-1}$ if necessary, one can 
arrange so that the pairs are the first ones 
$\{\alpha^{n}\beta^{m}\, ,\, \alpha^{-n}\beta^{-m}\}$ for 
$\forall (n, m) \in \mathbf Z^{2} - \{(0,0)\}$. 

\begin{claim} \label{claim:inf} There are real numbers 
$1 < C_{1} < C_{2}$ 
such that $\vert \mathcal T_{0} \vert = \infty$. Here 
$$\mathcal T_{0}\, :=\, \{(n, m)\, \in\, \mathbf Z^{2} - \{(0,0)\}\, \vert 
\, C_{1} < \alpha^{n}\beta^{m} < C_{2}\,\}\,\, .$$ 
\end{claim}

\begin{proof} Since $\alpha$, $\beta$, $C_{1}$, $C_{2} > 1$, the inequalty in 
$\mathcal T_{0}$ is equivalent to

$$(*)\,\, \text{log}\, C_{1} - n\,\text{log}\, \alpha < m\,\text{log}\,\beta < 
\text{log}\, C_{2} - 
n\,\text{log}\, \alpha \,\, .$$ 

Choose $1 < C_{1} < C_{2}$ such that $\beta < C_{2}/C_{1}$, i.e. 
$0 < \text{log}\, \beta < \text{log}\, C_{2} - \text{log}\, C_{1}$. 
For each $n \in \mathbf Z_{<0}$, the both sides of (*) 
are then positive, of distance $> \text{log}\, \beta$. Thus, for each such an 
$n$, there is at least one 
$m \in \mathbf Z_{>0}$ satisfying (*). This implies the result. 
\end{proof} 

We can now complete the proof of Proposition (4.4). Let us choose 
$C_{1}$ and $C_{2}$ as in Claim (4.5). Let $(n, m) \in \mathcal T_{0}$. 
Then the roots of $\varphi^{(n,m)}(t)$ are: 
$$\alpha^{n}\beta^{m}\,\, \alpha^{-n}\beta^{-m}\,\, 
\gamma_{i}\,\, \overline{\gamma_{i}}\,\, .$$ 
Here $1 \leq i \leq (\text{rank}\, T(M) - 2)/2$ and $\gamma_{i}$ 
are on the unit circle. One has  
$$C_{1} < \alpha^{n}\beta^{m} < C_{2}\,\, \text{and}\,\, 
1/C_{2} < \alpha^{-n}\beta^{-m} < 1/C_{1}\,\, .$$ 
Thus $\text{tr}\, \varphi^{(n,m)}(t)$ are bounded by 
the constants (independent on $(n, m) \in \mathcal T_{0}$). Hence, 
by Proposition (3.3), there is a set $\mathcal S$ of polynomials such that 
$\vert \mathcal S \vert < \infty$ and 
$$\varphi^{(n,m)}(t) \, \in \, \mathcal S 
\,\, \text{for}\,\, \forall (n, m) \, \in \, \mathcal T_{0}\,\, .$$
By $\vert \mathcal T_{0} \vert = \infty$ and $\vert \mathcal S \vert < \infty$, there is $\Psi(t) \in \mathcal S$ such that 
$\vert \mathcal T_{1} \vert = \infty$. Here   
$$\mathcal T_{1} := \{(n, m)\, \in \, \mathcal T_{0}\, \vert \, \varphi^{(n, m)}(t) = \Psi(t)\, \}\,\, .$$ 
The polynomial $\Psi(t)$ has at most finitely many roots. On the other hand, 
$\chi(f^{n}g^{m})$ is a root of $\varphi^{(n,m)}(t)$. Thus there is a root 
$\delta$ of 
$\Psi(t)$ 
such that $\vert \mathcal T_{2} \vert = \infty$. Here   
$$\mathcal T_{2} := \{(n, m)\, \in\, \mathcal T_{1}\, \vert\, \chi(f^{n}g^{m}) = \delta\}\,\, .$$ 
Choose one $(n_{0}, m_{0}) \in \mathcal T_{2}$.  
Then for each 
$(n, m) \in \mathcal T_{2}$, one has  
$$\chi(f^{n}g^{m}) = \chi(f^{n_{0}}g^{m_{0}})\,\, , \text{i.e.}\,\, 
\chi(f^{n -n_{0}}g^{m-m_{0}}) = 1\,\, .$$ 
Then $F^{n-n_{0}}G^{m-m_{0}} = id$ by $A \simeq A'$. 
Since $F$ and $G$ form a part of free basis of $A$, it follows that $n= n_{0}$ and 
$m = m_{0}$, i.e. $\mathcal T_{2} = \{(n_{0}, m_{0})\}$, 
a contradiction to $\vert \mathcal T_{2} \vert = \infty$.

Hence $r = 1$ if $r > 0$. This completes the proof of Proposition (4.4) 
and therefore the proof of the main part of Theorem (1.5)(1).

It remains to show the existence of $M$ with $r =1$ and $0$. By the main part of Theorem (1.5)(1), 
the pair $(M_{m}, f_{m})$ ($m \geq 2$) in 
Example (2.6) satisfies 
$r = 1$ 
as well as other necessary requirements. In addition, the hyperk\"ahler 
manifold $D_{m}$ in Corollary (1.7) satisfies $r = 0$ as well as other necessary requirements. 
$2$-dimensional example with $r = 0$ will be given in Section 6. 
\end{proof} 
\begin{remark}\label{remark:holy}
Here the Dirichlet unit theorem, one of two most fundamental theorems in algebraic number 
theory, played a crucial role in the proof. It may be interesting to seek a relation 
hidden between hyperk\"ahler manifolds and the other fundamental theorem, the finiteness of class numbers 
(if any). See [HLOY] for one of such relations in dimension $2$. 
\end{remark} 

\section{Automorphism group of a hyperk\"ahler manifold with parabolic N\'eron-Severi lattice}
In this section we shall prove the main part of Theorem (1.5)(2). We freely use some basic properties 
of almost abelian groups. They are in Section 9. Optimality of the estimate in dimension 
$2$ will be shown in Section 7. 
Throughout this section $M$ is a hyperk\"ahler manifold with parabolic $NS(M)$.

Set $H := {\rm Im}\,(r_{NS} : {\rm Bir}\,(M) \longrightarrow {\rm O}\,(NS(M)))$. 
By Corollary (2.7) and Proposition (9.3), we may show that $H$ is almost abelian of rank at most $\rho(X) - 1$. 
However, this follows from the next purely lattice theoretical:

\begin{proposition} \label{proposition:para} 
Let $L$ be a parabolic lattice of rank $r$. Let $N < {\rm O}\,(L)$. Then 
$N$ is an almost abelian group of rank at most $r-1$. 
\end{proposition}

\begin{proof} Let $v$ be the unique (up to sign) primitive 
totally isotoropic element of $L$. Then $g(v) = \pm v$ for each $g \in N$, 
and those $g$ with $g(v) = v$ form a subgroup of $N$ of index at most $2$. 
So, by Proposition (9.3), we may assume that $g(v) = v$ for all $g \in N$. 

The bilinear form of $L$ descends to the bilinear 
form on  $\overline{L} := L/\mathbf Z v$ and makes $\overline{L}$ 
a negative definite lattice of rank $r -1$. Our isometry 
$N$ also descends to 
the isometry of $\overline{L}$, say $g \mapsto \overline{g}$.

Set 
$$N^{(0)}\,\,  :=\,\, \text{Ker}\,\, 
(N \longrightarrow {\rm O}(\overline{L})\,\,;\,\,  g\, 
\mapsto\, \overline{g})\,\, .$$

By the negative definteness of $\overline{L}$, we have 
$[N : N^{(0)}] < \infty$. 
So, again by Proposition (9.3), it now suffices to show that 
$N^{(0)}$ is almost abelian of rank at most $r -1$. 
 
Let $\langle v, u_{i} (1 \leq i \leq r -1)\rangle$ be an 
integral basis of $L$. Let $g \in N^{(0)}$. Then we have 
$$g(v) = v\,\, ,\,\, g(u_{i}) = u_{i} + \alpha_{i}(g)v\,\, .$$ 

Here $\alpha_{i}(g)$ ($1 \leq i \leq r-1$) are integers uniquely determined 
by $g$. It is easy to see that the following map $\varphi$ is a group homomorphism:

$$\varphi : N \longrightarrow \mathbf Z^{r-1}\,\, ;\,\, g \mapsto 
(\alpha_{i}(g))_{i=1}^{r-1}\,\, .$$

This $\varphi$ is clearly injective by the form of $g$. Thus the result follows. 
\end{proof}
\begin{remark}\label{remark:geom} 
It might be more desirable to give a more geometric argument. It is somehow possible in dimension $2$. 
Indeed, if $S$ is a K3 surface with parabolic $NS(S)$, then the algebraic dimension $a(S) = 1$ and 
$S$ admits the unique elliptic fibration $a : S \longrightarrow \mathbf P^{1}$ (the algebraic reduction). This $f$ is 
${\rm Aut}\,(S)$-stable. Moreover, as it is easily seen, $f$ admits at least three singular 
fibers (cf. also [VZ]). Thus the group $K := {\rm Ker}\, ({\rm Aut}\, (S) \longrightarrow {\rm Aut}\, (\mathbf P^{1}))$ 
is of finite index in ${\rm Aut}\,(S)$. Since $f$ has no horizontal curve by $a(S) = 1$, 
the group $K$ acts on each fiber as translations. Thus, $K$ can be embedded into the Mordell-Weil group 
of the Jacobian fibration $j : J \longrightarrow \mathbf P^{1}$ of $f$. The surface $J$ is an algebraic K3 surface. 
Therefore $K$ is an abelian group of rank 
at most $18$ (cf. [Sh]). However, it seems difficult 
to generalize this argument in higher dimensional case, because almost nothing is known about the algebraic 
dimension and the structure of the algebraic reduction maps of $M$.
\end{remark}

\par \vskip 1pc
\section{Automorphism group of a McMullen's K3 surface}
In this section we shall prove Corollary (1.6) and the existence part in dimension $2$ 
of Theorem (1.5)(1). 
Let 
$(S, f)$ be a McMullen's pair.
Then $r = 1$ by the main part of Theorem (1.5)(1). First, we shall show that $\text{Aut}\,(S) 
\simeq \mathbf Z$. 
Again by the main part of Theorem (1.5)(1), this follows from the next Proposition which we learned 
from Professor JongHae Keum: 
\begin{proposition} \label{proposition:keum} {\rm [Keum]} Let 
$S$ be a K3 surface with $\rho(S) = 0$. Let $g \in \text{Aut}\,(S)$. Then 
$g = id$ iff $\chi(g)$ is a root of unity. 
\end{proposition}

\begin{proof} Only if part is trivial.  
Assume that $\chi(g)$ is a root of unity. Then $\chi(g) = 1$ by Theorem 
(2.4)(6) applied for $\langle g \rangle$. Thus 
$g^{*} \vert T(S) = id$ by Theorem (2.4)(2) 
and hence  $g^{*} = id$ on $H^{2}(S, \mathbf Z)$ by 
$\rho(S) = 0$. 
Then one has $g = id$ by the global Torelli theorem for K3 surfaces 
(see for instance [BPV]). 
\end{proof} 

The following proposition and Proposition (6.1) 
complete the proof of the rest: 

\begin{proposition} \label{proposition:generic} There are at most countably 
many K3 
surfaces 
with $\rho(S) = 0$ and having an automorphism $f$ such that $\chi(f)$ is not 
a root of unity. 
In particular, in the 
($20$-dimensional) period domain, 
K3 surfaces such that $\rho(S) = 0$ and ${\rm Aut}\,(S) = \{id\}$ form the 
complement of the 
union of countably many rational hyperplanes and at most countably 
many points.  
\end{proposition}
 
\begin{proof} By the surjectivity of the period mapping 
and the Lefschetz $(1,1)$-Theorem, K3 surfaces with $\rho (S) = 0$ 
form the complement of the unioin of the countably many rational 
hyperplane in the period 
domain. Thus last statement follows from the first statement and 
Proposition (6.1).

Let us show the first statement. The following argument is similar to [Mc]. 
Let $S$ be a K3 surface with $\rho(S) = 0$ and with $f \in \text{Aut}\, (S)$ 
such that $\chi(f)$ is not a root of unity. For each such 
$S$, let us choose 
a marking $\iota_{S} : T(S) = H^{2}(S, \mathbf Z) \simeq \Lambda$, where $\Lambda$ 
is the K3 lattice (cf. Section 7). The characteristic 
polynomial of $\iota_{S} \circ f^{*} \circ \iota_{S}^{-1}$  
must be a Salem 
polynomial by Theorem (3.4). There are at most countably 
many Salem polynomials 
of degree $22$. So, there are at most countably many element 
$g \in \text{O}(\Lambda)$ whose 
characteristic 
polynomials 
are Salem polynomials. By Corollary (3.11), each eigenvalue of such 
$g$ is 
of multiplicity one and the corresponding eigenspaces 
are one-dimensional. Thus, there are only countably many eigenspaces 
of all such $g$. On the other hand, the 
possible $1$-dimensional 
subspaces $\iota_{S}(\mathbf C \sigma_{S}) \subset \Lambda \otimes \mathbf C$ 
must be 
an eigenspace of one of such $g$. Thus, the subspaces 
$\iota_{S}(\mathbf C \sigma_{S})$ are at most countably many 
as well. 
Then, by the global Torelli Theorem, there are at most countably 
many such $S$. 
\end{proof} 
\par \vskip 1pc
\section{Non-projective K3 surface of large automorphism group} 

In this section, we shall complete the proof of Theorem (1.5)(2) by showing the following:

\begin{theorem} \label{theorem:large} 
There is a K3 surface $S$ such that $a(S) = 1$, $\rho(S) = 19$ and such that 
${\rm Aut}\,(S)$ is an almost abelin group of rank $18 (= 19 -1)$. 
In particular, the estimate in Theorem (1.5)(2) is optimal in dimension 
$2$. 
\end{theorem} 

It is well-known that $a(S) = 1$ iff $NS(S)$ is parabolic, 
and $\rho(S) \le 19$ if $a(S) \le 1$. On the other hand, by the main part of Theorem (1.5), 
${\rm Aut}\, (S)$ is almost abelian of rank at most 
${\rm max}\,(\rho(S) - 1, 1)$ if $a(S) \le 1$. So, the K3 surface $S$ in Theorem (7.1) is 
also a non-projective K3 
surface having maximal possible automorphism group (up to finite group factors).

In what follows, we shall show Theorem (7.1) dividing into several steps. 
Unfortunately, our proof is based on the surjectivity of the period mapping 
and the global Torelli Theorem so that our K3 surface and the group action 
are not so "visible".

Let 

$$\Lambda := \Lambda_{\text{K3}} := U^{\oplus 3} \oplus E_{8}(-1)^{\oplus 2}$$ 

be the K3 lattice. Let $i \in \{0, 1, 2\}$ and $\langle e_{i}, f_{i}\rangle$ be the integral basis 
of the $(i+1)-$th $U$ in $\Lambda$ s.t. $(e_{i}^{2}) = (f_{i}^{2}) = 0$, 
$(e_{i}, f_{i}) = 0$, and $\langle v_{ij} \rangle_{j=1}^{8}$ be an integral 
basis of the $i$-th $E_{8}(-1)$ ($i=1,2$) in $\Lambda$, which forms the Dynkin diagram of type 
$E_{8}$. (We need not specify the position of each $v_{ij}$ in the Dynkin 
diagram.)

Choose "sufficiently large" mutually different primes numbers $p$, $q$, $p_{i}$, $q_{i}$ ($1 \le i \le 8$). (The term 
"sufficiently large" will be 
clear in the proof of the next Lemma.)  We then define the sublattices $\overline{N}$, 
$N$ and $L$ of $\Lambda$ by:

$$\overline{N} := \mathbf Z \langle e_1 - pf_1 \rangle \oplus 
\mathbf Z \langle e_2 - qf_2 \rangle \oplus 
\mathbf Z \langle e_1 - p_{1j}v_{1j} \rangle_{j=1}^{8} \oplus 
\mathbf Z \langle e_2 - q_{2j}v_{2j} \rangle_{j=1}^{8}\,\, ,$$
$$N := \mathbf Z \langle e_0 \rangle \oplus \overline{N} \,\, ,$$ 
$$L := N + \mathbf Z\langle f_{0} \rangle = U \oplus \overline{N}\,\, .$$  
\begin{lemma} \label{lemma:ltice} 
$N$ and $L$ satisfy:
\begin{list}{}{
\setlength{\leftmargin}{10pt}
\setlength{\labelwidth}{6pt}
}
\item[(1)] $N$ is parabolic of rank $19$, and $L$ is a hyperbolic lattice 
of rank $20$.
\item[(2)] $N$ and $L$ are primitive in $\Lambda$.
\item[(3)] $N$ does not represent $-2$, i.e. there is no $x \in N$ with 
$(x^2) = -2$.
\end{list}
\end{lemma} 
\begin{proof} Note that $\overline{N}$ is elliptic of rank $18$ and that 
$(e_{0}^{2}) = (e_{0} \overline{N}) = 0$. 
Thus, $N$ is parabolic of rank $19$. (This does not depend on the choice of 
$p$, $q$, $p_{i}$, $q_{i}$ whenever they are positive integers.) Proof for 
$L$ is almost the same and we omit it.

Next we shall show that $N$ is primitive in $\Lambda$. (Proof for $L$ 
is the same and we omit it.) Here, we use the fact that 
$p$, $q$, $p_{i}$, $q_{i}$ are 
mutually different prime numbers (but not yet use the assumption they are "sufficiently large"). 
Suppose that there are rational numbers 
$x_0$, $x_1$, $x_2$, $x_{1i}$, $x_{2i}$ such that 
$$x_{0}e_{0} + x_{1}(e_1 - pf_1) + x_{2}(e_2 - qf_2) + 
\sum_{j=1}^{8} x_{1j}(e_1 - p_{1j}v_{1j}) + 
\sum_{j=1}^{8} x_{2j}(e_2 - p_{2j}v_{2j}) \in \Lambda\,\, .$$

We need to show that $x_{0}, x_{1}, x_{2}, x_{1j}, x_{2j}$ are all integers.

Calculate that
$$x_{0}e_{0} + x_{1}(e_1 - pf_1) + x_{2}(e_2 - qf_2) + 
\sum_{j=1}^{8} x_{1j}(e_1 - p_{1j}v_{1j}) + 
\sum_{j=1}^{8} x_{2j}(e_2 - p_{2j}v_{2j})$$
$$= x_{0}e_{0} + (x_{1} + \sum_{j=1}^{8}x_{1j})e_{1} - x_{1}pf_{1} + 
(x_{2} + \sum_{j=1}^{8}x_{2j})e_{2} - x_{2}qf_{2} - 
\sum_{j=1}^{8}x_{1j}p_{j}v_{1j} - \sum_{j=1}^{8}x_{2j}q_{j}v_{2j}\,\, .$$

Since $e_{0}, e_{1}, f_{1}, e_{2}, f_{2}, v_{1j}, v_{2j}$ form a part of free 
basis of $\Lambda$, the coefficients of the second line are all integers:

$$x_{0}\,\, ,\,\, x_{1}p\,\, ,\,\, x_{2}q\,\, ,\,\, x_{1j}p_{j}\,\, ,\,\, 
x_{2j}q_{j}\,\, \in \,\, \mathbf Z\,\, ;$$ 
$$x_{1} + \sum_{j=1}^{8}x_{1j}\,\, ,\,\,  x_{2} + \sum_{j=1}^{8}x_{2j}\,\, 
\in \,\, \mathbf Z\,\, .$$

Thus $x_{0} \in \mathbf Z$ and one can write $x_{1} = c_{1}/p$, $x_{2} = c_{2}/q$, $x_{1j} = c_{1j}/p_{j}$, 
$x_{2j} = c_{2j}/q_{j}$, where $c_{*}$ and $c_{**}$ are integers. By substituting these into the second two quantities, 
one obtains 
$$\frac{c_{1}}{p} + \sum_{j=1}^{8} \frac{c_{1j}}{p_{j}}\,\, ,\,\, 
\frac{c_{2}}{q} + \sum_{j=1}^{8} \frac{c_{2j}}{q_{j}}\,\, \in\,\, \mathbf Z\,\, .$$

By clearing the denominator of the first quantity, one has 
$$pp_{1} \cdots p_{8} \vert (c_{1}p_{1}\cdots p_{8} + 
\sum_{j=1}^{8} c_{1j}pp_{1} \cdots p_{j-1}p_{j+1} \cdots p_{8})\,\, .$$ 
Since $p$ divides the second sum, we have $p \vert c_{1}p_{1} \cdots p_{8}$, 
and therefore $p \vert c_{1}$. Similarly, $p_{j}$ must divide the term 
$c_{1j}pp_{1} \cdots p_{j-1}p_{j+1} \cdots p_{8}$, whence $p_{j} \vert c_{1j}$. 
For the same reason, one has also $q \vert c_{2}$ and $q_{j} \vert c_{2j}$. 
Hence the coefficients $x_{0}$, $x_{1}$, $x_{2}$, $x_{1j}$, $x_{2j}$ are all integers.

Let us show the assertion (3). Here we use the fact that $p$, $q$, $p_{i}$, $q_{i}$ are "sufficiently large". 
Since $E_{8}(-1)$ is negative definite, 
there are only finitely many 
elements $r$ with $(r^{2}) = -2$. 
(As wellknown, there are exactly $240$ such elements, but we do not need 
this precise numbers.) Then, there is no $(-2)$-element of the form 
$\sum_{j=1}^{8}x_{ij}n_{j}v_{ij}$ ($x_{ij} \in \mathbf Z$) if $n_{j}$ 
($1 \le j \le 8$) are sufficiently large, say $n_{j} > C$ for $\forall j$. 
Choose $8$ prime numbers $p_{j}$, $q_{j}$ so that they are larger than such 
$C$.

Let $v\in N$ and write $v$ under the integral basis as:
$$v =  x_{0}e_{0} + x_{1}(e_1 - pf_1) + x_{2}(e_2 - qf_2) + 
\sum_{j=1}^{8} x_{1j}(e_1 - p_{j}v_{1j}) + 
\sum_{j=1}^{8} x_{2j}(e_2 - q_{j}v_{2j})\,\, .$$

Then, by using $(e_{i}, v_{ij}) = (e_{i}^{2}) = 0$, one calculates:
$$(v^{2}) = -2px_{1}^{2} -2qx_{2}^{2} + (\sum_{j=1}^{8} x_{1j}p_{j}v_{1j})^{2} + (\sum_{j=1}^{8} x_{2j}q_{j}v_{2j})^{2}\,\, .$$

Here each of the four summands in $(v^{2})$ is even and $0$ or 
strictly less than $-2$ 
by our choice of $p$, $q$, $p_{j}$, $q_{j}$. Thus 
$(v^{2}) \not= -2$. This completes the proof. 

\end{proof} 

As usual, by a marked K3 surface, we mean a pair $(S, \iota)$ of a K3 surface 
$S$ and an isometry $\iota : H^{2}(S, \mathbf Z) \simeq \Lambda$. 

\begin{proposition} \label{proposition:exist} 
There is a marked K3 surface $(S, \iota)$ such that $\iota(NS(S)) = N$. 
\end{proposition} 

\begin{proof} 
Put 
$$T := N_{\Lambda}^{\perp} = \mathbf Z e_0 \oplus \overline{N}_{U^{\oplus 2} \oplus E_8(-1)^{\oplus 2}}^{\perp} = 
\mathbf Z e_{0} \oplus \overline{T}\,\, .$$ 
Then $\overline{T}$ is a positive definite lattice of rank $2$. Choose an integral basis $\langle u_{1}, u_{2} \rangle$ 
of $\overline{T}$. Put $(u_{1}, u_{1}) = 2a$, $(u_{2}, u_{2}) = 2c$, $(u_{1}, u_{2}) = b$, and $A := 4ac - b^{2}$. 
Then $a > 0$ and $A > 0$. Consider the element $\sigma$ of $T_{\mathbf C}$ defined by 
$$\sigma = \sqrt{2}e_{0} + \frac{-b + \sqrt{A}i}{2a}u_{1} + u_{2}\,\, i = \sqrt{-1}\,\, .$$ 
It is easy to check that 
$$(\sigma, \sigma) = 0\,\, , (\sigma, \overline{\sigma}) = \frac{A}{a} > 0\,\, .$$ 

Thus, by the sujectivity of the period mapping of K3 surfaces (see eg. [BPV]),  there is a marked K3 surface 
$(S, \iota)$ such that $\iota(\sigma_{S}) = \sigma$. Note that $N = T_{\Lambda}^{\perp}$. Then 
the equality $\iota(NS(S)) = N$ follows from the next Lemma. 
\end{proof} 
\begin{lemma} \label{lemma:trans} 
Under the identification of $H^{2}(S, \mathbf Z)$ with $\Lambda$ by $\iota$, the lattice $T$ 
is the transcendental lattice of $S$.  
\end{lemma} 
\begin{proof} Let $M$ be the minimal primitive sublattice of $\Lambda$ such that $\sigma \in M_{\mathbf C}$. 
It is clear that $M \subset T$. We need to show that $T \subset M$. 
Since $\sigma - \overline{\sigma} \in M_{\mathbf C}$, we have $(\sqrt{A}i/a)u_{1} \in M_{\mathbf C}$. Thus $u_{1} \in M$ 
(by the primitivity of $M$), and therefore $\sqrt{2}e_{0} + u_{2} \in M_{\mathbf C}$ by $\sigma \in M_{\mathbf C}$. 
Since $e_{0}, u_{2} \in \Lambda$, this implies $\sqrt{2}e_{0} + u_{2} \in M_{\mathbf Q(\sqrt{2})}$, and therefore 
$-\sqrt{2}e_{0} + u_{2} \in M_{\mathbf Q(\sqrt{2})}$ by the Galois theory. Thus $e_{0}, u_{2} \in M$ again by the primitivity 
of $M$. Hence $T \subset M$. 
\end{proof} 
The next theorem will complete the proof. 
\begin{theorem} \label{theorem:dis} 
Let $S$ be a K3 surface in Proposition (7.3). Then $a(S) = 1$, $\rho(S) = 19$ and ${\rm Aut}\, (S)$ 
is an almost abelian group of rank $18$. 
\end{theorem} 
In what follows, we identify $H^{2}(S, \mathbf Z)$ with $\Lambda$ by the marking $\iota$ and use the notation in 
Proposition (7.3) and its proof. 
Then $NS(S) = N$ and $T(S) = T$. Thus $a(S) = 1$ and $\rho(S) = 19$. Let us consider ${\rm Aut}(S)$. 
Let $L$ be the hyperbolic lattice defined just before Lemma (7.2). We have $L_{\Lambda}^{\perp} = \overline{T}$. 
\begin{lemma} \label{lemma:iso} For each $i$ ($1 \le i \le 18$), there is an isometry $\varphi_{i}$ of $L$ of the following form:
$$e_{0} \mapsto e_{0}\,\, ,\,\, f_{0} \mapsto f_{0} + \gamma_{i}e_{0} + \sum_{k=1}^{18}c_{ik}w_{ik}\,\, ;$$
$$w_{i} \mapsto w_{i} + m w_{i}\,\, ;\,\, w_{j} \mapsto w_{j}\,\, (j \not= i)\,\, .$$
\end{lemma} 
\begin{proof} Put $Q := ((w_{k}, w_{l}))_{1 \le k, l \le 18}$. This is a negative definite integral matrix. 
Put $m = {\rm det}\, Q$. Then $\varphi_{i}$ induces an isometry of $N$ and satsfies 
$(\varphi_{i}(f_{0}), \varphi_{i}(e_{0}) = (f_{0}, e_{0})$ (regardless with values $\gamma_{i}, c_{ik}$). 
Now it suffices to find $\gamma_{i}, c_{ik} \in \mathbf Z$ such that 
$$(\varphi_{i}(f_{0}), \varphi_{i}(w_{k})) = (f_{0}, w_{k})\,\, (1 \le k \le 18)\,\, ,\,\, 
(\varphi_{i}(f_{0})^{2}) = (f_{0}^{2})\,\, .$$
Since $m = \text{det}\, Q$, the first $18$ equalities are euivalent to that $(c_{ik})_{k=1}^{18}$ is the (minus) the $i$-th low 
of the adjoint matrix $\tilde{Q}$ of $Q$. Thus $c_{ik} \in \mathbf Z$. The last equality is equivalent to 
$$2\gamma_{i} + ((\sum_{k=1}^{18}c_{ik}w_{k})^{2}) = 0\,\, .$$ 
Since $\Lambda$ is even, we have $\gamma_{i} \in \mathbf Z$. Now $\varphi_{i}$ is an isometry of $L$ 
for these $c_{ik}$ and $\gamma_{i}$. 
\end{proof} 
\begin{lemma} \label{lemma:ext} Let $\varphi_{i}$ be an isometry of $L$ in Lemma (7.6). Then, for each $i$ ($1 \le i \le 18$), 
there is a positive 
integer $k_{i}$ such that the isometry $(\varphi_{i}^{k_{i}}, id_{\overline{T}}) \in {\rm O}(L) \times {\rm O}(\overline{T})$ 
extends to an isometry $\Phi_{i}$ of $\Lambda$. 
\end{lemma} 
\begin{proof} Let $L^{*}$ is the dual (over-) lattice of $L$. We may choose $k_{i}$ so that the induced action of 
$\varphi_{i}^{k_{i}}$ on $L^{*}/L$ is identity. So, we may put $k_{i} = \vert L^{*}/L \vert$. 
\end{proof}
\begin{proposition} \label{proposition:au} Let $\Phi_{i}$ be an isometry of $\Lambda$ in Lemma (7.7). 
Then, for each $i$ ($1 \le i \le 18$), there is 
an automorphism $F_{i} \in {\rm Aut}\, (S)$ such that $F_{i}^{*} = \Phi_{i}$.  
\end{proposition} 
\begin{proof} By the global Torelli theorem for K3 surfaces (see eg. [BPV]), we may show that 
$\Phi_{i}$ is an effective Hodge isometry of $\Lambda = H^{2}(S, \mathbf Z)$. By the construction, $\Phi_{i}(T) = T$ and $\Phi_{i} \vert T = id_{T}$. 
In particular $\Phi_{i}(\sigma_{S}) = \sigma_{S}$. Thus $\Phi_{i}$ is an Hodge isometry of $H^{2}(S, \mathbf Z)$. 
Recall that $N = NS(S)$ does not represent $(-2)$-element. Thus, 
$S$ contains no smooth rational curve. It follows that the K\"ahler cone $\mathcal K(S)$ of $S$ coincides with the positive cone 
$\mathcal C$ of 
$H^{1,1}(S, \mathbf R)$. By replacing $e_{0}$ by $-e_{0}$ (if necessary), we may assume that $e_{0}$ is in the boundary 
of $\mathcal C$. Since $\Phi_{i}(e_{0}) = e_{0}$, we have then $\Phi_{i}(\mathcal C) = \mathcal C$. Hence 
$\Phi_{i}(\mathcal K(S)) = \mathcal K(S)$ and $\Phi_{i}$ is also effective. 
\end{proof} 
The next proposition will complete the proof of Theorem (7.5) and hence that of Theorem (7.1):
\begin{proposition} \label{proposition:al} Let $F_{i}$ be the automorphism of $S$ in Proposition (7.8). Then: 
\begin{list}{}{
\setlength{\leftmargin}{10pt}
\setlength{\labelwidth}{6pt}
}
\item[(1)] Let 
$G := \langle F_{i} \rangle_{i=1}^{18}$ be a subgroup of ${\rm Aut}\,(S)$ generated by $F_{i}$'s. 
Then $G$ is a free abelian group of rank $18$.
\item[(2)] ${\rm Aut}\, (S)$ is an almost abelian group of rank $18$.
\end{list} 
\end{proposition} 
\begin{proof} Let $g \in G$. Then, we have $g(L) = L$ and more precisely 
$g^{*}(e_{0}) = e_{0}$, $g^{*}(w_{i}) = w_{i} + m_{i}(g)w_{i}$ and 
$g^{*}(f_{0}) = f_{0} + \gamma(g)e_{0} + \sum c_{k}(g)w_{k}$ for some integers $m_{i}(g)$ $\gamma(g)$ and $c_{k}(g)$. 
This is because the corresponding equality holds for the generators $F_{i}^{*}$ by $F_{i}^{*}\vert L = \varphi_{i}$ (see 
Lemma (7.6)). 
Then, we have a group homomorphism $\alpha : A \longrightarrow \mathbf Z^{18}$ defined by $g \mapsto (m_{i}(g))_{i=1}^{18}$. 
Note also that $g^{*} \vert \overline{T} = id_{\overline{T}}$, because so are $F_{i}^{*}$. 
Using these informations, one can show that $\alpha(g) = 0$ iff $g^{*} \vert \Lambda = id$ in the same manner as in Proposition 
(5.1). However, $g^{*} \vert \Lambda = id$ iff $g = id$ by the global Torelli theorem. Thus $\alpha$ is injective. 
Combining this with the fact that $\alpha(F_{i}) = (\delta_{ik}m)_{k=1}^{18} \in \mathbf Z^{18}$, where $m$ is a positive integer 
in Lemma (7.6), we obtain the assertion (1).

Let us show (2). We freely use the facts about almost abelian group in Section 9. 
By the main part of Theorem (1.5)(2), ${\rm Aut}\,(S)$ is almost abelian of rank $r \le 18$, 
i.e. there is a normal subgroup $H$ of ${\rm Aut}\,(S)$ such that 
$[G:H] < \infty$ and $H$ fits in with the exact sequence 
$$1 \longrightarrow N \longrightarrow H \longrightarrow \mathbf Z^{r} \longrightarrow 1\,\, ,$$ 
where $N$ is a finite subgroup of $H$. We need to show that $r = 18$. The exact sequence above induces the following exact 
sequence 
$$1 \longrightarrow N \cap G \longrightarrow H \cap G \longrightarrow \mathbf Z^{r}\,\, .$$ 
Since $G$ is a free abelian group, we have $N \cap G = \{0\}$. Thus $H \cap G$ is a subgroup of $\mathbf Z^{r}$. 
On the other hand, since $H$ is of finite index in ${\rm Aut}\, (S)$ and $F_{i} \in {\rm Aut}\,(S)$, there is a positive 
integer $n$ such that 
$$F_{1}^{n}\,\, ,\,\, F_{2}^{n}\,\, ,\,\, \cdots \,\, ,\,\, F_{18}^{n} \in H\,\, ,$$ 
and thus $H \cap G$ contains a subgroup isomorphic to $\mathbf Z^{18}$. 
Hence $18 \le r$, and therefore $r = 18$. 
\end{proof}  
\par \vskip 1pc

\section{New counterexamples of Kodaira's problem}

In this section, we shall prove Theorem (1.9). 

Let $(S, f)$ be a McMullen's pair. Put $Y := S \times S$. Let $D \subset Y$ 
be the diagonal and $G \subset Y$ be the graph of $f^{-1}$. 
we have $D \not= G$. Let $x : X \longrightarrow Y$ 
be a composite of blow up along smooth centers lying over 
$D \cap G$ such that the proper transforms $D'$ and $G'$ (in $X$)
of $D$ and $G$ are smooth and disjoint. Let 
$z : Z := Z_{4} \longrightarrow X$ be 
the blow 
up along $D' \cup G'$. This $Z$ is a $4$-dimensional 
non-projective compact K\"ahler manifold being 
bimeromorphic to $S \times S$. In particular, $\pi_{1}(Z) = \{1\}$ 
and $a(Z) = 0$. Here $a(Z)$ is the algebraic dimension of $Z$.

We first show the following: 

\begin{proposition} \label{proposition:rigid} 
$Z$ is a counterexample of Kodaira's problem in dimension $4$. 
\end{proposition}

\begin{proof} Let $\pi : (\mathcal Z \supset Z) \longrightarrow 
(\mathcal B \ni 0)$ be a small deformation of $Z \simeq \mathcal Z_{0}$ 
over a positive dimensional analytic space  
$\mathcal B$. For the purpose, by taking a (reduction, normalization and) 
resolution of 
$\mathcal B$ [Hi] and considering the pullback family, one may assume that 
$\mathcal B$ is smooth. 
In addition, since the problem is local, one may also assume 
that $0 \in \mathcal B$ 
is a (smooth) germ 
and one can freely shrink $0 \in \mathcal B$ 
whenever it will be convenient. Let $E(D')$ and $E(G')$ be the exceptinal 
divisors of $z$ lying over $D'$ and $G'$. Then, by Kodaira [Ko], 
there are smooth subfamilies $\mathcal E(D') \subset \mathcal Z$ 
and $\mathcal E(G') \subset \mathcal Z$ over $\mathcal B$ such that 
$\mathcal E(D')_{0} = E(D')$ 
and $\mathcal E(G')_{0} = E(G')$. Then, by a result of Fujiki and Nakano 
[FN], $\mathcal E(D')$ and 
$\mathcal E(G')$ 
can be contracted to smooth subfamilies $\mathcal D' \subset \mathcal X$ and 
$\mathcal G' \subset \mathcal X$ over $\mathcal B$. Then 
$\mathcal X_{0} = X$, $\mathcal D'_{0} = D'$ and $\mathcal G_{0}' = G'$. 
Let $E \subset X$ be the exceptional divisor of the last blow up 
$X \longrightarrow X'$ in $x$. 
Then, for the same reason, $E$ is extended over $\mathcal B$ and then 
is simultaneously contracted to a smooth family $\mathcal X'$, say 
$x' : 
\mathcal X \longrightarrow \mathcal X'$. Under this contraction, 
$\mathcal D' \longrightarrow \mathcal D"$ is also smooth 
blow-down (possibly isomorphism), 
because so is $D' \longrightarrow D"$. (Here $D"$ is the proper transform 
of $D$ on $X'$.) The same holds for 
$\mathcal G' \longrightarrow \mathcal G"$.  
Repeating this process, one finally obtains a smooth proper family 
$$\varphi : \mathcal Y \longrightarrow \mathcal B$$ and its two smooth 
proper subfamilies over $\mathcal B$: 
$$\mathcal D, \mathcal G \subset \mathcal Y$$ 
such that 
$$\mathcal Y_{0} = Y\, ,\, \mathcal D_{0} = D\, ,\, \mathcal G_{0} = G\, .$$ 
Note that smooth proper defomation of a K3 surface is again a K3 surface. 
Thus $\mathcal D_{b}$ and $\mathcal G_{b}$ ($b \in \mathcal B$) 
are both K3 surfaces, because so are $D$ and $G$. Recall 
that $Y = S \times T$ (where for convenience, we write the second factor 
by $T$, which is actually $S$). Since 
$h^{1}(\mathcal O_{S}) = h^{1}(\mathcal O_{T}) = 0$, the fibered manifold 
structures 
$p_{1} : S \times T \longrightarrow S$ and 
$p_{2} : S \times T \longrightarrow T$ are both stable 
under deformations by Kodaira [ibid]. More precisely, 
there are smooth proper families 
$$p : \mathcal S \longrightarrow 
\mathcal B\,\, ,\,\, q : \mathcal T \longrightarrow \mathcal B\,\, $$ 
such that $S \simeq \mathcal S_{0}$ and $T \simeq \mathcal T_{0}$,  
and smooth proper families over $\mathcal B$  
$$r : \mathcal Y \longrightarrow \mathcal S
\,\, ,\,\, s : \mathcal Y \longrightarrow 
\mathcal T\,\, ,$$ 
such that $p \circ r = q \circ s = \varphi$. For the same reason as before, 
$\mathcal S_{b}$ and $\mathcal T_{b}$ are all K3 surfaces. 

Since the fiber $T$ of $r$ over $0 \in \mathcal B$ meets each of 
$D$ and $G$ at 
one point transeversally, the same is also true for the fiber of $r$ over 
$b \in \mathcal B$ 
(near $0$). The same is true for $s$. 
 Thus the natural morphisms 
$$r\vert \mathcal D_{b} : 
\mathcal D_{b} \longrightarrow \mathcal S_{b}\,\, , \,\, 
s\vert \mathcal D_{b} : \mathcal D_{b} \longrightarrow 
\mathcal T_{b}\,\, $$ 
$$r\vert \mathcal G_{b} : \mathcal G_{b} \longrightarrow 
\mathcal S_{b}\,\, , \,\, 
s\vert \mathcal G_{b} : \mathcal G_{b} \longrightarrow 
\mathcal T_{b}\,\, $$ 
are all bimeromorphic morphisms among K3 surfaces, and hence are all 
isomorphisms.

Then the morphism $\tilde{f} : \mathcal S \longrightarrow \mathcal S$ 
defined by
$$\mathcal S \mapright{(r \vert \mathcal D)^{-1}} \mathcal D 
\mapright{s \vert \mathcal D} \mathcal T 
\mapright{(s \vert \mathcal G)^{-1}} \mathcal G \mapright{r 
\vert \mathcal G} \mathcal S$$ 
is an automorphism of $\mathcal S$ over $\mathcal B$ such that 
$\tilde{f}_{0} = f$. This $\tilde{f}$ naturally acts on the local system 
$$R^{2}p_{*} \mathbf Z_{\mathcal S}$$ 
as an isometry $\tilde{f}^{*}$, and also acts on the Hodge filtration
$$p_{*}\Omega_{\mathcal S/\mathcal B}^{2} \subset 
R^{2}p_{*} \mathbf Z_{\mathcal S} \otimes \mathcal O_{\mathcal B}\,\, .$$ 
Let $\tilde{\sigma}$ is a non-zero 
section of $p_{*}\Omega_{\mathcal S/\mathcal B}^{2}$. Then 
$\tilde{f}^{*}\tilde{\sigma} = \tau \tilde{\sigma}$. Here $\tau$ 
is an element of $\Gamma(\mathcal B, \mathcal O_{\mathcal B}^{\times})$. 
One has  
$\tau (b) = \chi(\tilde{f}_{b})$ and $\tau(0) = 
\chi(\tilde{f}_{0}) = \chi(f)$. 
On the other hand, since 
$R^{2}p_{*} \mathbf Z_{\mathcal S}$ is a constant system, the matrix 
representation of $\tilde{f}$ with respect to a flat basis of 
$R^{2}p_{*}\mathbf Z_{\mathcal S}$, is constant over 
$\mathcal B$. Thus, the eigenvalues of 
$\tilde{f}_{b} \in O(H^{2}(\mathcal S_{b}, \mathbf Z))$ are also constant 
over $\mathcal B$. In particular, 
$\chi(\tilde{f}_{b}) = \tau(b)$ are constant and therefore equal to $\chi(f)$. 
Since $(S, f)$ is McMullen's pair, $\chi(f)$ is not a root of unity. Thus, 
$\mathcal D_{b} \simeq \mathcal S_{b}$ is not projective by Theorem (2.4). 
Hence $\mathcal Y_{b}\, (\supset \mathcal D_{b})$ is not projective, either. 
On the other hand, $\mathcal Y_{b}$ 
is K\"ahler, for being K\"ahler is stable under small proper 
deformations [KS]. 
Thus $a(\mathcal Y_{b}) < \text{dim}\, \mathcal Y_{b}$ 
by a result of Moishezon [Mo]. 
Since $\mathcal Z_{b}$ is bimeromorphic to $\mathcal Y_{b}$, one has then 
$a(\mathcal Z_{b}) < \text{dim}\, \mathcal Z_{b}$. 
Thus $\mathcal Z_{b}$ is not projective, and $Z$ is not algebraically 
approximated. 
\end{proof} 
\begin{remark} \label{remark:strong} Recall that projecive K3 surfaces 
are dense in any non-trivial deformation of a K3 surface [Fu2] (See also [GHJ, Proposition 26.6]). On the other 
hand, $\mathcal S_{b}$ are all non-projective as observed above. 
Thus, the family $\mathcal S \longrightarrow \mathcal B$ 
is a constant family and $\mathcal S_{b}$ are all isomorphic to $S$. 
Then the pairs 
$(\mathcal S_{b}, \tilde{f}_{b})$ are also isomorphic to $(S, f)$ by the 
discreteness of $\text{Aut}\, (S)$ or by 
$\chi(\tilde{f}_{b}) = \chi(f)$ and 
$\text{Aut}\, S \simeq \mathbf Z$ (Corollary (1.6)). 
Since $\mathcal Z_{b}$ is uniquely recovered from the pair 
$(\mathcal S_{b}, \tilde{f}_{b})$, it follows that $\mathcal Z_{b}$ are all 
isomorphic to $Z$, i.e. $Z$ is analytically rigid. 
\end{remark}

Let $Z = Z_{4}$ be as in Proposition (8.1) and $F_{m}$ be a smooth 
hypersurface of degree $m + 2$ in $\mathbf P^{m+1}$. 
Let $d \geq 6$ and set $Z_{d} := Z \times F_{d-4}$. Then $Z_{d}$ is bimeromorphic 
to a 
$d$-dimensional simply-connected K\"ahler manifold with trivial canonical 
bundle. This $Z_{d}$ is not rigid, 
because the second factor $F_{d-4}$, 
which is the algebraic reduction 
of $Z_{d}$, actualy moves in $\mathbf P^{d-3}$. 
We also set $Z_{5} := Z \times \mathbf P^{1}$. Then $\pi_{1}(Z_{d}) = \{1\}$
for $\forall d \geq 4$.

Now the following proposition completes the proof of Theorem (1.9). 

\begin{proposition} \label{proposition:simplerigid} 
Let $d \geq 5$. Then $Z_{d}$ is a 
counterexample of Kodaira's problem in dimension $d$.
\end{proposition}

\begin{proof} Let $\pi : (\mathcal Z_{d} \supset Z_{d}) \longrightarrow 
(\mathcal B \ni0)$ be any small deformation 
of $Z_{d}$. As in Proposition (8.1), one may assume that 
$0 \in \mathcal B$ is a smooth germ. Since $Z_{d} = Z \times F_{d-4}$ and 
$h^{1}(\mathcal O_{Z}) = 0$, 
again by the stability theorem of Kodaira, we have a 
family $p : (\mathcal F_{d-4} \supset F_{d-4})  \longrightarrow 
(\mathcal B \ni 0)$ and a family of fibered manifolds 
$r : \mathcal Z \longrightarrow \mathcal F_{d-4}$ 
over $\mathcal B$ such that $\pi = p \circ r$. 
Let $\tilde{0} \in \mathcal F_{d-4}$ such that $p(\tilde{0}) = 0$. 
Then $r^{-1}(\tilde{0}) \simeq Z$ by the construction. Thus, by 
Proposition (8.1), $r^{-1}(c)$ is non-projective for 
$\forall c \in \mathcal U$. Here $\mathcal U$ is 
an open subset of $\mathcal F_{d-4}$ such that $\tilde{0} \in \mathcal U$. 
Thus the fibers 
$$\pi^{-1}(p(c))\, =\, (\mathcal Z_{d})_{p(c)}\, (\supset r^{-1}(c))\,$$ 
are non-projective for $\forall c \in \mathcal U$. Since $p$ is a smooth 
surjective 
morphism, $p(\mathcal U)$ is an open subset of $\mathcal B$ 
with $0 \in p(\mathcal U)$. This completes the proof. 
\end{proof}
\vskip 1pc 
\section{Almost abelian group (Appendix)}
This section is a sort of appendix, in order to make statements and arguments in 
the main part (Sections 1-8) clear. We shall give a precise definition of almost abelian group (we employed) 
and remark a few easy properties of almost abelian groups (we used in the main part). 

\begin{definition} \label{definition:alb} A group $G$ is {\it almost abelian group} 
(resp. {\it almost 
abelian of finite rank}, say $r$) 
if there are a normal subgroup $G^{(0)}$ of $G$ of finite index, 
a finite group $K$ and an abelian group $A$ (resp. a non-negative integer $r$)  which fit in 
the exact sequence 
$$1 \longrightarrow K \longrightarrow G^{(0)} \longrightarrow A \longrightarrow 0\,\,  $$ 
(resp. 
$$1 \longrightarrow K \longrightarrow G^{(0)} \longrightarrow 
\mathbf Z^{r} \longrightarrow 0\,\, ).$$  
\end{definition} 

It is clear that if $G$ is almost abelian, then so are subgroups $H < G$ and quotient groups $G/N$. 

First, we show the well-definedness of the rank of an almost abelian group. 

\begin{lemma} \label{lemma:rank} 
Let $G$ be an almost abelian group of finite rank. Then 
the rank of $G$ is uniquely determined by $G$.
\end{lemma} 

\begin{proof} Let $G$, $G^{0}$ and $r$ be as in the definition (9.1). 
Let $G^{(1)}$ be another normal subgroup of $G$ of finite index which 
fits in an exact sequence
$$1 \longrightarrow N \longrightarrow G^{(1)} \longrightarrow 
\mathbf Z^{s} \longrightarrow 0\,\,  .$$ 
Here $N$ is a finite group. We need to show that $r = s$. 
Set $G^{(2)} = G^{(0)} \cap G^{(1)}$. 
From the standard exact sequence
$$1 \longrightarrow G^{(0)} \longrightarrow G \longrightarrow G/G^{(0)} 
\longrightarrow 1\,\, ,$$ 
we have an exact sequence 
$$1  \longrightarrow G^{(2)} \longrightarrow G^{(1)} 
\longrightarrow G/G^{(0)}\,\, .$$ 
Thus, $G^{(2)}$ is a normal subgroup of $G^{(1)}$ of finite index. 
From the first exact sequence, we have an exact sequence 
$$1 \longrightarrow N \cap G^{(2)} \longrightarrow G^{(2)} 
\longrightarrow \mathbf Z^{s}\,\, .$$ 
Since $N$ is a finite group and $[G^{(1)}:G^{(2)}] < \infty$, 
the image of the last map 
is also of finite index in $\mathbf Z^{s}$. Thus, it is isomorphic to 
$\mathbf Z^{s}$ 
(as abstract abelian groups) and we have an exact sequence 
$$1 \longrightarrow N^{(2)} \longrightarrow G^{(2)} \longrightarrow 
\mathbf Z^{s} \longrightarrow 0\,\,  ,$$ 
where $N^{(2)}$ is a finite subgroup of $G^{(2)}$.  
Similarly, we have 
an exact sequence
$$1 \longrightarrow K^{(2)} \longrightarrow G^{(2)} \longrightarrow 
\mathbf Z^{r} \longrightarrow 0\,\, ,$$ 
where $K^{(2)}$ is a finite subgroup of $G^{(2)}$. 
Taking ${\rm Hom}(*, \mathbf Z)$ as groups, we obtain exact sequences 
$$0 \longrightarrow {\rm Hom}(\mathbf Z^{r}, \mathbf Z) 
\longrightarrow {\rm Hom}(G^{(2)}, \mathbf Z) 
\longrightarrow {\rm Hom}(K^{(2)}, \mathbf Z)\,\, = 0\,\, ,$$ 
$$0 \longrightarrow {\rm Hom}(\mathbf Z^{s}, \mathbf Z) 
\longrightarrow {\rm Hom}(G^{(2)}, \mathbf Z) 
\longrightarrow {\rm Hom}(N^{(2)}, \mathbf Z)\,\, = \,\, 0\,\, .$$ 
For the last two equalities, we used the fact that $K^{(2)}$ 
and $N^{(2)}$ are finite. Thus,  
$$\mathbf Z^{r} \simeq  {\rm Hom}(G^{(2)}, \mathbf Z) \simeq 
\mathbf Z^{s}\,\, .$$ 
Hence $r = s$ by the fundamental theorem of finitely generated abelian groups. 
\end{proof} 

The next proposition has been frequently used in the proof of Theorem (1.5):
\begin{proposition} \label{proposition:kerim} 
Let $G$ be a group.
\begin{list}{}{
\setlength{\leftmargin}{10pt}
\setlength{\labelwidth}{6pt}
}
\item[(1)] Let
$$1 \longrightarrow H \longrightarrow G \longrightarrow Q \longrightarrow 1$$ 
be an exact sequence such that $Q$ is a finite group. Then, $H$ 
is an almost abelian group 
if and only if  so is $G$. Moreover the ranks of $H$ and $G$ are the same. 
\item[(2)]Let 
$$1 \longrightarrow K \longrightarrow G \longrightarrow H 
\longrightarrow 1$$ 
be an exact sequence such that $K$ is a finite group. Then, $H$ 
is an almost abelian 
group if and only if so is $G$. Moreover the ranks of $H$ and $G$ 
are the same. 
\end{list}
\end{proposition} 

\begin{proof} We shall show (1). The proof of (2) is similar. 
Assume that $G$ is almost abelian. Take a normal subgroup $G^{(0)}$ as 
in the definition (9.1). 
Set $H^{(0)} := G^{(0)} \cap H$. Then, as in the proof of Lemma 
(9.2), $H^{(0)}$ makes $H$ almost abelian 
of the same rank as $G$. This shows "if part". Assume that $H$ is almost abelian 
and take a normal subgroup $H^{(0)}$ of $H$ as in the definition (9.1).  
Then 
$$G^{(0)} := \cap_{[g] \in G/H}\,\, g^{-1}H^{(0)}g$$ 
is a normal subgroup of $G$ of finite index 
(cf. [Su, Chapter 1]). 
As before, $G^{(0)}$ makes $G$ almost 
abelian of the same rank of $H$. This shows "only if part".  
\end{proof} 

%\vfill\eject
%%%%%%%%%%%%%% reference %%%%%%%%%%%%%%%%%

\end{document}